\documentclass[11pt,twoside]{amsart}
\usepackage{amsfonts}
\usepackage{epsfig,graphics}
\usepackage{amssymb}
\usepackage{amscd}

\newtheorem*{theor}{Th\'eor\`eme}
\newtheorem{theoreme}{Th\'eor\`eme}[section]

\newtheorem{proposition}[theoreme]{Proposition}

\newtheorem{lemme}[theoreme]{Lemme}

\newtheorem{question}[theoreme]{Question}

\newtheorem{remarque}[theoreme]{Remarque}

\newtheorem{corollaire}[theoreme]{Corollaire}

\newcommand{\diag}{\text{diag}}

\newcommand{\dd}{\text{d}}

\newcommand{\OO}{\text{O}}
\newcommand{\CO}{\text{CO}}
\newcommand{\PO}{\text{PO}}

\newcommand{\Ein}{\text{Ein}}

\newcommand{\Conf}{\text{Conf}}

\newcommand{\Ad}{\text{Ad }}

\newcommand{\vecconf}{\chi^{co}(M)}
\newcommand{\Aff}{\text{Aff}}
\newcommand{\End}{\text{End}}

\newcommand{\BN}{{\bf N}} 
\newcommand{\BR}{{\bf R}}
\newcommand{\BC}{{\bf C}}

\newcommand{\RP}{{\bf RP}}
\newcommand{\PL}{{\bf PL}}

\newcommand{\lien}{{\mathfrak{n}}}

\newcommand{\lieg}{{\mathfrak{g}}}

\newcommand{\liep}{{\mathfrak{p}}}

\newcommand{\lier}{{\mathfrak{r}}}
\newcommand{\oo}{{\mathfrak{o}}}

\begin{document}

\pagenumbering{arabic}
\title[Formes normales pour les champs conformes]{Formes normales pour les
  champs conformes pseudo-riemanniens}
\author{Charles Frances et Karin Melnick}
\date{\today}
\maketitle

\begin{abstract}
We establish normal forms for conformal vector fields on pseudo-Riemannian manifolds in
the neighborhood of a singularity.  For real-analytic Lorentzian manifolds, we show that
 the vector field is analytically linearizable or the manifold is conformally
flat.  In either case, the vector field is locally conjugate to a normal form on a model
space. For smooth metrics of general signature, we obtain the analogous result under the
additional assumption that the differential of the flow at the fixed point is bounded.
\end{abstract}


\section{Introduction}

Cet article \'etudie l'allure locale des champs de vecteurs conformes,
c'est-\`a-dire les champs de vecteurs lisses $X$ sur une vari\'et\'e
pseudo-riemannienne $(M,g)$ qui satisfont ${\mathcal L}_Xg=\sigma g$,
pour une certaine fonction $\sigma \in C^{\infty}(M)$. Ce sujet a d\'ej\`a fait l'objet d'une litt\'erature abondante,   notamment dans le cadre des  champs de
vecteurs conformes sur les espaces-temps, c'est-\`a-dire les vari\'et\'es
lorentziennes (voir en particulier \cite{beig}, \cite{capocci1},
\cite{capocci2}, \cite{kuehnel},  \cite{kuehnel2},  \cite{kuehnel4}, et \cite{steller}).  

\`A titre d'exemple, commen\c{c}ons par d\'ecrire une famille int\'eressante de champs conformes, qui va jouer un r\^ole fondamental dans tout l'article.  Pour tout couple d'entiers naturels $(p,q)$, $1\leq p \leq q$, le projectivis\'e du c\^one de lumi\`ere d'une forme quadratique de signature $(p+1,q+1)$ est naturellement muni d'une structure pseudo-riemannienne conforme de signature $(p,q)$.  Ce projectivis\'e, muni de cette structure conforme naturelle est un espace compact  appell\'e {\it l'univers d' Einstein} de signature $(p,q)$. On le note ${\Ein}^{p,q}$. Le groupe $\PO(p+1,q+1)$ agit transitivement sur ${\Ein}^{p,q}$ par transformations conformes.  Soit $o$ un point fix\'e sur ${\Ein}^{p,q}$. Le stabilisateur de $o$ dans  $\PO(p+1,q+1)$ est un sous-groupe parabolique  que nous notons $P$. Ainsi,  du point de vue conforme, ${\Ein}^{p,q}$ est simplement l'espace homog\`ene $\PO(p+1,q+1)/P$.  Maintenant, tout groupe \`a un param\`etre de $P$ donne lieu \`a un champ de vecteurs conforme sur ${\Ein}^{p,q}$, admettant une singularit\'e $o$.  De tels champs conformes seront appel\'es {\it champs de M\"obius}.  Nous verrons que d\'ej\`a dans le cadre lorentzien, il existe de nombreuses  formes locales diff\'erentes pour les champs de M\"obius au voisinage d'une de leurs singularit\'es.

Soit maintenant $(M,g)$ une vari\'et\'e pseudo-riemannienne de signature $(p,q)$, avec $p+q \geq 3$, et $X$ un champ de vecteurs conforme sur $M$ s'annulant en un point $x_0$. L'existence d'une connexion de  Cartan canonique associ\'ee \`a la structure conforme $[g]$  permet d'\'etudier le champ $X$ gr\^ace \`a la connaissance de son jet d'ordre 2 en $x_0$.  Cette propri\'et\'e a d\'ej\`a \'et\'e utilis\'ee avec succ\`es, par exemple dans \cite{aleks}, ou \cite{nagano} pour l'\'etude des champs de vecteurs projectifs.  Dans 
\cite{nilpotent}, puis dans \cite{frances-riemannien}, nous avons syst\'ematis\'e cet usage de la connexion de Cartan pour l'\'etude des champs conformes, et nous avons expliqu\'e
comment associer \`a $X$ un champ de M\"obius $X_h$ admettant une singularit\'e en $o$, que l'on appelle le {\it champ d'holonomie de $X$ en $x_0$}.  Dans ces pr\'ec\'edents travaux, nous avons commenc\'e \`a \'etablir   un dictionnaire entre les propri\'et\'es locales de $X_h$ au voisinage de $o$, et celles de $X$ au voisinage de $x_0$.  Ceci nous a conduit \`a formuler la question suivante~:

\begin{question}
\label{quest.conjugaison}
En dimension $n \geq 3$, un champ de vecteurs conforme pseudo-riemannien $X$, est-il toujours localement conjugu\'e, au voisinage d'une singularit\'e $x_0$, \`a son champ d'holonomie? 
\end{question}

Bien entendu, une r\'eponse positive \`a cette question r\`eglerait de mani\`ere tr\`es satisfaisante le probl\`eme de trouver des formes normales pour les champs conformes au voisinage d'une  singularit\'e.   On serait en effet ramen\'e \`a l'\'etude des champs de M\"obius, pour lesquels on peut faire des calculs explicites.
 
Rappelons  que dans le cadre riemannien, la r\'eponse  \`a la Question \ref{quest.conjugaison} est affirmative.  Cela d\'ecoule du th\'eor\`eme ci-dessous, qui peut s'obtenir \`a partir des travaux de D. Alekseevskii (voir \cite{aleks}, Section $2$).  Une  preuve d\'etaill\'ee se trouve dans \cite{frances-riemannien},  Th\'eor\`eme 1.2.

\begin{theor} 
{\it Soit $(M,g)$ une vari\'et\'e riemannienne, lisse, de dimension $n \geq 3$,  munie d'un champ de vecteurs conforme $X$.  On suppose que $X$ s'annule en $x_0 \in M$. Alors ou bien $X$ est lin\'earisable, et complet  sur un voisinage de $x_0$, ou bien il existe un voisinage de $x_0$ qui est conform\'ement plat.
Dans tous les cas, le champ $X$ est $C^{\infty}$-conjugu\'e au voisinage de $x_0$ \`a son champ d'holonomie $X_h$.}
\end{theor}




\subsection{\'Enonc\'es des r\'esultats}

Le but du pr\'esent article est de r\'epondre, au moins partiellement, \`a la Question \ref{quest.conjugaison} dans le cadre des signatures autres que riemanniennes.  On obtient en particulier un r\'esultat complet pour les structures conformes lorentziennes  analytiques~:

\begin{theoreme}
\label{thm.lorentz-analytique}
Soit $(M,g)$  une vari\'et\'e lorentzienne, analytique, de dimension $n \geq 3$.  Soit $X$ un champ de vecteurs conforme analytique sur $M$, admettant une singularit\'e $x_0$. Alors :

\begin{enumerate}
\item{Ou bien $X$ est analytiquement lin\'earisable au voisinage de $x_0$.}
\item{Ou bien $(M,g)$ est conform\'ement plate.}
\end{enumerate}
Dans les deux cas, le champ $X$ est analytiquement conjugu\'e au voisinage de $x_0$ \`a son champ d'holonomie $X_h$.
\end{theoreme}

Dans certains cas, on peut  obtenir un r\'esultat de m\^eme nature, en signature $(p,q)$  g\'en\'erale, et toujours pour des structures analytiques.  Avant de l'\'enoncer, pr\'ecisons que si $X$ est un champ de vecteurs conformes ayant une singularit\'e $x_0$, alors la diff\'erentielle du flot local de $X$ en $x_0$, not\'ee $D_{x_0}\phi_X^t$, est d\'efinie pour tout $t \in {\bf R}$, et d\'etermine, \`a conjugaison pr\`es, un sous-groupe \`a un param\`etre de transformations conformes lin\'eaires de ${\text CO}(p,q)$.

\begin{theoreme}
\label{thm.semi-simple-analytique}
Soit $(M,g)$ une vari\'et\'e analytique, pseudo-riemannienne, de dimension $n \geq 3$, et $X$ un champ de vecteurs conforme sur $M$,  admettant une singularit\'e en $x_0$. On suppose que $\{D_{x_0}\phi_X^t \}_{t \in \BR}$  est un groupe \`a un param\`etre semi-simple sur $\BC$.  Alors~:
\begin{enumerate}
\item{Ou bien $X$ est analytiquement lin\'earisable au voisinage de $x_0$.}
\item{Ou bien $(M,g)$ est conform\'ement plate.}
\end{enumerate}
Dans les deux cas, $X$ est analytiquement conjugu\'e au voisinage de $x_0$ \`a son champ d'holonomie $X_h$.
\end{theoreme} 

Sans faire d'hypoth\`ese d'analyticit\'e, nos m\'ethodes permettent encore d'obte\-nir des informations sur l'allure locale de certains champs conformes.   Un champs conforme sur $(M,g)$ est \emph{essentiel} si son flot local ne pr\'eserve aucune m\'etrique de la classe conforme $[g]$. Il est dit {\it inessentiel} sinon. 
\begin{theoreme} 
\label{thm.translations}
Soit $(M,g)$ une vari\'et\'e lisse, pseudo-riemannienne, de dimension $n \geq 3$, et $X$ un champ de vecteurs conforme sur $M$,  admettant une singularit\'e en $x_0$.  On note $\{ \phi_X^t \}$ le flot local engendr\'e par $X$ sur $M$, et on suppose que $\{ D_{x_0}\phi_X^t \}_{t \in \BR}$  est un groupe \`a un param\`etre relativement compact de $\CO (T_{x_0} M)$.  Alors :
\begin{enumerate}
\item{Ou bien il existe un voisinage $U$ de $x_0$ sur lequel $X$  est complet  et engendre un flot relativement compact dans $\Conf(U)$.  Dans ce cas $X$ est lin\'earisable en $x_0$, et est inessentiel sur $U$.}
\item{Si l'on n'est pas dans le premier cas, $X$ est essentiel sur tout voisinage de $x_0$.  Dans ce cas, $(M,g)$ poss\`ede  un ouvert conform\'ement plat et contenant  $x_0$ dans son adh\'erence.}
\end{enumerate}
\end{theoreme}

Ce r\'esultat est une \'etape essentielle dans les preuves des Th\'eor\`emes \ref{thm.lorentz-analytique} et \ref{thm.semi-simple-analytique}.
Nous verrons que dans le cas $(2)$ du th\'eor\`eme, on a de plus une
description dynamique relativement pr\'ecise du flot $\{ \phi_X^t \}$ au voisinage de la singularit\'e $x_0$, 
ainsi
qu'une description g\'eom\'etrique de l'ouvert o\`u $[g]$ est
conform\'ement plate (voir les \'enonc\'es des Th\'eor\`emes
\ref{thm.translation-partiel} et \ref{thm.translations-partiel-lumiere}).

Il est naturel de se demander s'il ne serait pas possible d'am\'eliorer les conclusions dans le cas  $(2)$ du th\'eor\`eme en montrant qu'un voisinage de la singularit\'e $x_0$ doit \^etre conform\'ement plat.   Il n'en est rien, par exemple lorsque la signature est lorentzienne, comme le montrent les exemples de la Section 6 de \cite{frances-causal}.  Toutefois, si l'on fait une hypoth\`ese de compacit\'e sur $M$, on peut raisonablement esp\'erer des conclusions plus fortes. En particulier~:

\begin{question}
Soient $(M,g)$ et $X$ comme dans le Th\'eor\`eme \ref{thm.translations}.  On suppose que $M$ est compacte, et que le flot global   $\{ \phi^t_X \}$ n'est pas relativement compact dans $\Conf(M)$.  La vari\'et\'e  $(M,g)$ est-elle forc\'ement  conform\'ement plate?
\end{question}

Dans cette direction, le second auteur
 et Andreas \v Cap ont r\'ecemment \'etudi\'e des flots essentiels sur des
g\'eom\'etries paraboliques  g\'en\'erales dans \cite{cap.me.parabolictrans}. Ils ont \'etendu des
outils d\'evelopp\'es dans \cite{nagano} et  le pr\'esent article,
et les ont
 appliqu\'es \`a diverses g\'eom\'etries  pour obtenir des
 descriptions dynamiques au voisinage de points fixes et montrer, dans
certains cas,  l'annulation de la
courbure.

\subsection{R\'esum\'e de l'article.  }  
Les preuves des th\'eor\`emes de cet article sont bas\'ees sur l'id\'ee que certains comportements dynamiques locaux de transformations conformes ne peuvent advenir que sur des structures conform\'ement plates.  Par exemple, en signature riemannienne, si une transformation conforme est une contraction topologique sur un ouvert, on montre par un argument classique que cet ouvert est conform\'ement plat.  Toutefois, la situation se complique en signature quelconque~: nous donnerons, en derni\`ere section de l'article, un exemple de contraction lin\'eaire agissant conform\'ement pour une structure pseudo-riemannienne qui n'est pas conform\'ement plate.   Par ailleurs, \`a partir de la signature lorentzienne, les motifs dynamiques  locaux des transformations conformes sont plus nombreux et subtils que de simples contractions topologiques.  Il faut donc des outils pour comprendre pr\'ecis\'ement  la dynamique conforme locale.

Le lecteur trouvera en  Section \ref{sec.rappels-geometriques} quelques rappels sur les structures conformes pseudo-riemanniennes, et leur interpr\'etation comme g\'eom\'etries de Cartan.  
Puis, dans la Section \ref{sec.holonomie}, nous  introduisons notre outil dynamique principal~: l'holonomie associ\'ee \`a une suite de transformations conformes. Il s'agit d'une suite de transformations conformes du mod\`ele $\Ein^{p,q}$, qui permet d'\'etudier la dynamique des transformations conformes d'origine  (voir la Proposition \ref{prop.stabilite} et le Th\'eor\`eme \ref{thm.integrabilite-stable}).
Les preuves des Th\'eor\`emes \ref{thm.translations} et \ref{thm.semi-simple-analytique} sont donn\'ees dans la Section \ref{sec.preuve-translations}.   Celle du Th\'eor\`eme \ref{thm.lorentz-analytique} est l'objet de la Section \ref{sec.lorentz-analytique}.  Il s'agit de  la plus longue, car elle utilise tous les r\'esultats ant\'erieurs, ainsi qu'une connaissance assez pouss\'ee de la dynamique conforme dans l'espace mod\`ele lorentzien $\Ein^{1,n-1}$.  

\section{Rappels g\'eom\'etriques}
\label{sec.rappels-geometriques}

 
Cette section est consacr\'ee \`a une br\`eve description de la g\'eom\'etrie de l'espace mod\`ele $\Ein^{p,q}$, ainsi qu'\`a l'interpr\'etation des structures conformes en dimension $\geq 3$ en termes de g\'eom\'etrie de Cartan.  On note $n = p+q$ et on suppose que $p \leq q$.
 \subsection{L'espace mod\`ele : univers d'Einstein} 
Consid\'erons $\BR^{p+1,q+1}$,  l'espa\-ce $\BR^{p+q+2}$ muni de la forme quadratique :
$$ Q^{p+1,q+1}(x):=2x_0x_{p+q+1}+ \cdots +2x_px_{q+1}+\Sigma_{p+1}^q x_i^2$$
dont le c\^one isotrope est not\'e ${\mathcal N}^{p+1,q+1}$.  La restriction de $Q^{p+1,q+1}$ fournit une m\'etrique d\'eg\'en\'er\'ee sur ${\mathcal N}^{p+1,q+1} \setminus \{0 \}$, dont le noyau est de dimension $1$, et est tangent aux g\'en\'eratrices du c\^one.  Ainsi, le projectivis\'e ${\bf P}({\mathcal N}^{p+1,q+1} \setminus \{ 0 \})$ est une sous-vari\'et\'e lisse de $\RP^{p+q+1}$, naturellement munie d'une classe conforme de m\'etriques non d\'eg\'en\'er\'ees, de signature $(p,q)$.  
On appelle {\it univers d'Einstein} de signature $(p,q)$, et l'on note
$\Ein^{p,q}$, cette vari\'et\'e compacte ${\bf P}({\mathcal N}^{p+1,q+1}
\setminus \{  0  \})$ munie de la structure conforme ci-dessus.  Il s'agit
d'un espace remarquable au sein des structures conformes de signature $(p,q)$,
et nous renvoyons le lecteur \`a \cite{primer} et \cite[Chap. 4]{charlesthese} pour une \'etude d\'etaill\'ee dans le 
cadre lorentzien, ainsi qu'\`a  \cite[Sec. 3]{nilpotent} pour  le cas g\'en\'eral de la signature $(p,q)$.

Commen\c{c}ons par dire quelques mots sur le groupe des transformations
conformes de $\Ein^{p,q}$.  Si $\OO(p+1,q+1)$ d\'esigne le groupe des
applications lin\'eaires qui laissent $Q^{p+1,q+1}$ invariante, alors
l'action naturelle de $\PO(p+1,q+1)$ sur $\Ein^{p,q}$ pr\'eserve
clairement  la classe conforme de $\Ein^{p,q}$. Il s'av\`ere que
$\PO(p+1,q+1)$ est en fait  tout le groupe des transformations conformes
de $\Ein^{p,q}$; c'est le contenu du th\'eor\`eme de Liouville.  L'action de $\PO(p+1,q+1)$ est transitive sur $\Ein^{p,q}$.  On peut donc repr\'esenter $\Ein^{p,q}$ comme un espace homog\`ene $\PO(p+1,q+1)/P$, o\`u $P$ est un sous-groupe parabolique de $\PO(p+1,q+1)$, que nous supposerons \^etre par la suite le stabilisateur du point $o=[e_0]$.  

\subsubsection{L'alg\`ebre de Lie $\oo(p+1,q+1)$}
\label{sec.algebre}
L'alg\`ebre de Lie $\oo(p+1,q+1)$ est compos\'ee des matrices $X$ de taille
$(n+2) \times (n+2)$ qui satisfont l'identit\'e : 
$$ X^t J_{p+1,q+1} + J_{p+1,q+1} X = 0.$$
Ici,  $J_{p+1,q+1}$ est la matrice, dans la base $(e_0, \ldots ,e_{n+1})$, de la forme quadratique $Q^{p+1,q+1}$.

L'alg\`ebre $\oo(p+1,q+1)$ s'\'ecrit comme une somme $\lien^- \oplus \lier \oplus \lien^+$, o\`u :

$$ \lier =  \left\{ \left( \begin{array}{ccc}
a &  & 0 \\
  & M   &  \\
  &    & -a
\end{array} \right) \ :    
\qquad 
\begin{array}{c}
 a \in \BR  \\
  M \in \oo(p,q) \\
    \end{array} 
\right\}
$$

$$ \lien^+= \left\{ \left( \begin{array}{ccc}
0& -x^t.J_{p,q}   &  0\\
  & 0  & x \\
  &   & 0
\end{array} \right) \ :  
\qquad 
\begin{array}{c}
  x\in \BR^{p,q} 
\end{array}
\right\}
$$

$$ \lien^- = \left\{ \left( \begin{array}{ccc}
0&    &  \\
 x & 0  &  \\
 0 & -x^t.J_{p,q}  & 0
\end{array} \right) \ :  
\qquad 
\begin{array}{c}
  x\in \BR^{p,q} 
\end{array}
\right\}
$$

On va aussi consid\'erer la sous-alg\`ebre ${\mathfrak a} \subset \oo(p+1,q+1)$, constitu\'ee des matrices~:
$${\mathfrak a} = \left\{ \left( \begin{array}{ccccccc}
\alpha_1&    & & & & &  \\ $$
 & \ddots & & & & &  \\
 & & \alpha_{p+1} & & &  &\\
 & & & 0_{q-p} & & &  \\
 & & & & -\alpha_{p+1}& & \\
 & & & & & \ddots & \\
 & & & & & & -\alpha_{1} \\
 \end{array}
\right)  :  
\qquad 
\begin{array}{c}
  \alpha_1, \ldots ,\alpha_{p+1} \in \BR 
\end{array} \right\}.
$$

Le sous groupe ferm\'e de $P$ dont l'alg\`ebre de Lie est ${\mathfrak a}$ est not\'e $A$. On appelle ${\mathfrak a}^+$ le sous-ensemble de ${\mathfrak a}$ pour lequel $\alpha_1\geq \ldots \geq  \alpha_{p+1} \geq 0$, et $A^+:= e^{{\mathfrak a}^+}.$  

\subsubsection{Cartes conformes}
Deux cartes int\'eressantes nous seront utiles par la suite.  La premi\`ere est $j : \BR^{p,q} \to \Ein^{p,q}$ donn\'ee en coordonn\'ees projectives sur ${\bf RP}^{p+1,q+1}$ par :

$$ j : (x_1,\ldots,x_n) \to [-\frac{Q^{p,q}(x)}{2}:x_1: \ldots : x_n : 1 ]$$

Il s'agit d'un plongement conforme qui envoie l'espace de Minkowski $\BR^{p,q}$ sur un ouvert dense de $\Ein^{p,q}$.  Ainsi $\Ein^{p,q}$ est un espace localement conform\'ement plat.  Le compl\'ementaire de $j(\BR^{p,q})$ dans $\Ein^{p,q}$ est le c\^one de lumi\`ere de sommet $o$, c'est-\`a-dire l'ensemble des g\'eod\'esiques de type lumi\`ere passant par $o$.  L'ouvert $j(\BR^{p,q})$ est invariant par $P$, et l'application $j$ conjugue l'action affine de $\CO(p,q) \ltimes \BR^n$ sur $\BR^{p,q}$ \`a l'action de $P$ sur $j(\BR^{p,q})$; en particulier, $P \cong \CO(p,q) \ltimes \BR^n$.  Le groupe $A$ sera souvent vu, via ces cartes, comme un sous-groupe de transformations conformes lin\'eaires et diagonales  de $\BR^{p,q}$. Les \'el\'ements de $A^+$ sont des transformations  de $\CO(p,q)$ de la forme~:
$$ h=\diag(\lambda_1, \ldots, \lambda_n), \text{ avec } \lambda_i   \geq 1, \text{ pour } i=1,\ldots,n.$$

Notons que la carte $j$ ne contient pas le point $o$, ni son c\^one de lumi\`ere. Nous avons une seconde carte, dite \emph{carte \`a l'infini}:
 $$ j^o: (x_1,\ldots,x_n ) \mapsto [1 : x_1 : \cdots: x_n : -\frac{Q^{p,q}(x)}{2}].$$
L'application $j^o$ est un diff\'eomorphisme conforme de $\BR^{p,q}$ sur un ouvert de  $\Ein^{p,q}$ contenant $o$.

Notons que $j$ et $j^o$ ne suffisent pas \`a recouvrir $\Ein^{p,q}$  (sauf dans le cas de la sph\`ere, $p=0$), mais ce ne sera pas un probl\`eme pour la suite.


\subsection{Structures conformes et g\'eom\'etries de Cartan}
\label{sec.principe-equivalence}
Dans ce qui suit, nous allons utiliser intensivement le fait qu'en dimension
$n \geq 3$, les structures pseudo-riemanniennes conformes sont ce que l'on
appelle des {\it g\'eom\'etries de Cartan}.  Ceci est r\'esum\'e par le th\'eor\`eme suivant, appel\'e aussi {\it principe d'\'equivalence}, dont le lecteur pourra trouver une preuve dans \cite[Th. 4.2]{kobayashi} ou \cite[Chap. 7, Prop. 3.1]{sharpe} :

\begin{theoreme}[E. Cartan]
\label{thm.principe-equivalence}
Soit $(M,[g])$ une structure conforme pseudo-riemannienne de signature $(p,q)$, $p+q \geq 3$. Alors $(M,[g])$ d\'efinit de mani\`ere  unique un triplet $(M,B,\omega)$ o\`u :
\begin{enumerate}
\item{$B$ est un $P$-fibr\'e principal au-dessus de $M$.}
\item{La forme $\omega$ est une $1$-forme \`a valeurs dans ${\mathfrak o}(p+1,q+1)$, satisfaisant:
\begin{enumerate}
\item{ Pour chaque $b \in B$, $\omega_b : T_bB \to \lieg$ est un isomorphisme d'espaces vectoriels. }
\item{Pour tout $p \in P$, si $R_p$ d\'esigne l'action \`a droite de $p$ sur $B$,  $(R_p)^*\omega = (\Ad p^{-1}) \circ \omega$.}
\item{Pour tout $X \in \lieg$ et $b \in B$, $ \omega_b(\left. \frac{d}{dt} \right|_{t =0}R_{e^{tX}}.b)=X.$}
\end{enumerate}}
\end{enumerate}
R\'eciproquement, tout triplet $(M,B,\omega)$ comme ci-dessus d\'efinit une structure conforme de signature $(p,q)$ sur $(M,[g])$.
\end{theoreme}

Dans le cas du mod\`ele $\Ein^{p,q}= \PO(p+1,q+1)/P$, le $P$-fibr\'e principal
$B$ est le groupe de Lie $G= \PO(p+1,q+1)$, et la forme $\omega$ est la forme de
Maurer-Cartan sur $G$, qui sera not\'ee $\omega_G$ par la suite.  
On notera $\pi_G$ la projection $G \rightarrow \Ein^{p,q} = G/P$.  Le
th\'eor\`eme pr\'ec\'edent dit que toute structure conforme peut \^etre interpr\'et\'ee comme une version courbe de ce mod\`ele plat.



Dans la suite de l'article, le triplet $(M,B,\omega)$ donn\'e par le Th\'eor\`eme \ref{thm.principe-equivalence} sera appel\'e {\it fibr\'e normal de Cartan} associ\'e \`a la structure conforme $(M,[g])$.

\subsection{Application exponentielle}
\label{sec.exponentielle}

Soit $(M,[g])$ une structure conforme de signature $(p,q)$, $p+q \geq 3$, et
$(M,B,\omega)$ le fibr\'e normal de Cartan qui lui est associ\'e.  Tout choix
d'un vecteur $Z \in \lieg$ d\'efinit un champ de vecteurs $\hat Z$ sur $B$,
caract\'eris\'e par $\omega(\hat Z) \equiv Z$.  Ceci permet de d\'efinir une
application exponentielle sur un voisinage $\mathcal{W}$ de $B \times \{0\}$ dans $B \times \lieg$ : 
\begin{eqnarray*}
\exp & : & \mathcal{W} \rightarrow B \\ 
\exp(b,Z) & = & \phi_{\hat Z}^1.b
\end{eqnarray*}
o\`u $\phi_{\hat Z}^t$ d\'esigne le flot local associ\'e au champ de vecteurs $\hat Z$.

Si $\lambda : I \to B$ et $p : I \to P$ sont deux  courbes de classe $C^1$, et si l'on pose $\gamma(s)=\lambda(s).p(s)$, alors :
\begin{equation}
\label{equ.formule-Cartan}
\omega(\gamma^{\prime}(s))=(\Ad p(s))^{-1} \circ \omega(\lambda^{\prime}(s))+\omega_G(p^{\prime}(s)).
\end{equation}

Cette relation d\'ecoule du troisi\`eme point de (2), dans le Th\'eor\`eme \ref{thm.principe-equivalence} (voir \cite[Ch 5]{sharpe}).  Si l'on prend $p(s)$ constant \'egal \`a $p$, on voit que
\begin{eqnarray}
\label{equ.equivariance-exponentielle}
\exp(b.p, t (\Ad p)^{-1}.\xi) = \exp(b,t \xi).p.
\end{eqnarray}

Par ailleurs, la connexion de Cartan
$\omega$ conduit \`a la notion fondamentale de d\'eveloppement des courbes,
d'une structure conforme $(M,[g])$ dans le mod\`ele $\Ein^{p,q}$.  Soit $I$ un intervalle $[0, a]$.  Pour $b \in
B$ soit $\hat \alpha : I \to B$ une courbe de classe $C^1$, telle que $\hat
\alpha(0)=b$.  On appelle {\it d\'eveloppement de $\hat \alpha$} en $b$,
not\'e ${\mathcal D}_b(\hat \alpha)$ l'unique courbe $\hat \beta : I \to G$
satisfaisant  $\hat \beta (0)=1_G$ et
$\omega(\alpha^{\prime}) \equiv \omega_G(\beta^{\prime})$.  Soit maintenant $x \in M$, $b \in B$ au-dessus de $x$, et $\alpha : I \to M$ une courbe de classe $C^1$ telle que $\alpha(0)=x$, alors on d\'efinit son \emph{d\'eveloppement en $x$, relativement \`a $b$} par~:
$${\mathcal D}_x^b(\alpha) :=\pi_G \circ {\mathcal D}_b(\hat \alpha),$$

o\`u $\hat \alpha$ est un relev\'e de $\alpha$ dans $B$, tel que $\hat \alpha(0)=b$.
La d\'efinition de ${\mathcal D}_x^b(\alpha)$ est ind\'ependante du relev\'e $\hat \alpha$ choisi \cite[Prop 5.4.13]{sharpe}.

\subsubsection{G\'eod\'esiques conformes}
\label{sec.geodesiques-conformes}

Nous appellerons {\it g\'eod\'esique conforme para- m\'etr\'ee} toute courbe de
la forme $s \mapsto \pi(\exp(b,s\xi))$, o\`u $s$ d\'ecrit un intervalle $I = (-\delta,\delta)$, et $\xi$ appartient \`a  $(\Ad P).\lien^-$.  Le support g\'eom\'etrique d'une g\'eod\'esique $\alpha$ sera not\'e $[\alpha]$ dans ce qui suit.


Des expressions  matricielles donn\'ees en Section  \ref{sec.algebre}, on d\'eduit facilement que l'image par $j$ d'une droite de $\BR^{p,q}$ param\'etr\'ee affinement  est une g\'eod\'esi\-que conforme param\'etr\'ee de $\Ein^{p,q}$.  Toutes les g\'eod\'esiques conformes para- m\'etr\'ees de $\Ein^{p,q}$ d\'efinis sur tout $\BR$ sont obtenues en composant les g\'eod\'esiques ci-dessus avec des \'el\'ements de $\PO(p+1,q+1)$.   Nous voyons donc que le vecteur tangent $v$ aux points d'un segment g\'eod\'esique conforme est de  type constant: {\it type temps} si $\langle v , v \rangle < 0$ par rapport  \`a n'importe quelle m\'etrique de la classe conforme de $\Ein^{p,q}$, {\it type espace} si $\langle v, v \rangle> 0$, et {\it type lumi\`ere} si $\langle v, v \rangle = 0$.  Enfin, observons que les segments g\'eod\'esiques conformes de type temps, ou espace, issus de $o$, auxquels on a \^ot\'e le point $o$, sont les images par $j$ des demi-droites de  $\BR^{p,q}$.

Remarquons que toutes les  g\'eod\'esiques conformes compl\`etes de $\Ein^{p,q}$ se compactifient en des cercles lisses, par adjonction d'un point.  Dans le cas des g\'eod\'esiques de type lumi\`ere, ces compactifications sont les projections sur $\Ein^{p,q}$ des plans totalement isotropes de $\BR^{p+1,q+1}$.  On d\'esignera souvent ces objets comme \'etant les {\it g\'eod\'esiques de lumi\`ere} (non param\'etr\'ees) {\it de $\Ein^{p,q}$}.


\subsubsection{Normes et m\'etriques auxiliaires}
\label{sec.normes}
On munit $\BR^{p,q}$ du produit scalaire euclidien 
$$b(x,y):=x_1y_1+ \cdots + x_ny_n,$$

et l'on note $||x||$ la norme associ\'ee.  On transporte cette m\'etrique euclidienne  par l'application $j^o$, et l'on obtient une m\'etrique $\rho^o$ sur $j^o(\BR^{p,q})$.  Cette m\'etrique induit, par l'action simplement transitive de $N^- < G$ sur $j^o(\BR^{p,q})$, un produit scalaire sur $\lien^-$, et une norme que l'on note $|| \cdot ||_{\lien^-}$.  Dans la suite, on notera ${\mathcal B}(0,r)$ la boule de centre $0$ et de rayon $r$ dans $\lien^-$, relativement \`a la norme $|| \cdot ||_{\lien^-}$ et ${\mathcal S}(0,r)$ la sph\`ere correspondante.  On notera \'egalement $B(o,r):=e^{{\mathcal B}(0,r)}$ et $S(o,r):=e^{{\mathcal S}(0,r)}$.  Remarquons que $B(o,r)$ et $S(o,r)$ sont, respectivement, la boule et la sph\`ere de rayon $r$, centr\'ees en $o$,  pour la m\'etrique $\rho^o$.

\subsubsection{Transformations conformes et application exponentielle}
\label{sec.transformations}
Soit $(M,[g])$ une structure conforme pseudo-riemannienne, et $(M,B,\omega)$ le fibr\'e normal de Cartan.

Si $x \in M$ et $b \in B$ est dans la fibre de $x$, on a un isomorphisme naturel~:
$$ \iota_b: \lieg/\liep \to T_xM$$
que l'on d\'efinit par $\iota_b(\overline{\xi}):=D_b\pi(\omega_b^{-1}(\xi))$, o\`u $\xi$ est n'importe quel repr\'esentant de la classe $\overline{\xi} \in \lieg / \liep$.

Les relations d'\'equivariance de $\omega$ impliquent que pour tout $b \in B$, $p \in P$:
\begin{equation}
\label{eq.iota}
\iota_{b.p^{-1}}((\overline{\mbox{Ad}}\  p).\overline{\xi})=\iota_b(\overline{\xi}),
\end{equation}
o\`u l'on a not\'e $\overline{\mbox{Ad}}$ l'action adjointe sur $\lieg/\liep$.

Par ailleurs, toute transformation conforme locale $f$ d'un ouvert de $M$ dans $M$ se remonte en un automorphisme local du fibr\'e $B$, qui satisfait $f^*\omega=\omega$.  Nous utiliserons fr\'equemment les deux relations fondamentales suivantes~:

\begin{equation}
\label{eq.f.exponentielle}
f(\exp(b,\xi)).p^{-1}=\exp(f(b).p^{-1}, (\Ad p). \xi),
\end{equation}

ainsi que 
\begin{equation}
\label{eq.differentielle.f}
D_xf(\iota_b(\xi))=\iota_{f(b)}(\xi).
\end{equation}

\subsubsection{Ensembles convexes}
\label{sec.convexe}
Soit $\lambda^{p,q}$ un produit scalaire de signature $(p,q)$ sur $\lien^-$, qui soit invariant par l'action adjointe de $\OO(p,q) < P$.   Par construction m\^eme du fibr\'e normal de Cartan, pour tout $x \in M$ et $b \in B$ dans la fibre de $x$, on a $\iota_b^*([g_x])=[\lambda^{p,q}]$, avec $[g_x]$ la classe conforme sur $T_xM$ et $[\lambda^{p,q}]$ la classe des m\'etriques $\alpha \lambda^{p,q}$, $\alpha \in \BR_+^*$, sur $\lien^-$.

Soit ${\mathcal U}$ un ouvert convexe de $\lien^-$ contenant $0$.  Le c\^one isotrope ${\mathcal N}(\lien^-)$ de $\lambda^{p,q}$ divise $\mathcal{U}$ en deux ouverts : 
\begin{itemize}
\item l'ouvert ${\mathcal U}^+:= \{ \xi \in {\mathcal U} \ | \ \lambda^{p,q}(\xi,\xi) >0  \}$, qui est aussi connexe.\\
\item  l'ouvert ${\mathcal U}^-:= \{ \xi \in {\mathcal U} \ | \ \lambda^{p,q}(\xi,\xi) <0  \}$.  Cet ouvert est connexe si $p \geq 2$, vide si $p=0$, et poss\`ede deux composantes connexes si $p=1$.
\end{itemize}

Il est clair que ${\mathcal U}^+$ et ${\mathcal U}^-$ demeurent inchang\'es si l'on remplace $\lambda^{p,q}$ par $\alpha \lambda^{p,q}$ avec $\alpha \in \BR^*$.  

On dira qu'un ouvert $U \subset  M$ est {\it convexe} s'il est diff\'eomorphe \`a un ouvert convexe ${\mathcal U} \subset \lien^-$ contenant $0$ via l'application  $\pi \circ \exp(b, \cdot)$, avec $b \in B$.
Si $U=\exp(b,{\mathcal U})$ est un ouvert convexe, on notera $U^+:=\exp(b,{\mathcal U}^+)$ et $U^-:=\exp(b,{\mathcal U}^-)$.   Par exemple, pour $r$ petit, les voisinages $B(o,r)$ sont convexes, et on notera $B^\pm(o,r)=e^{{\mathcal B}(0,r)^\pm}$.

L'alg\`ebre des champs de vecteurs conforme sur $M$ sera not\'ee $\vecconf$ dans la suite.

\section{Holonomie}
\label{sec.holonomie}

Dans toute cette section, on d\'esigne par $(M,g)$ une vari\'et\'e pseudo-riemannien\-ne de signature $(p,q)$, $p+q \geq 3$.  Le fibr\'e normal de Cartan associ\'e \`a la structure conforme $(M,[g])$ est not\'e $(M,B,\omega)$.  On consid\`ere  $X \in \vecconf$ admettant une singularit\'e $x_0 \in M$.

\subsection{Champ d'holonomie et suite d'holonomie}
\label{sec.holonomie-champs}
 Par construction du fibr\'e normal de Cartan, toute transformation conforme locale $f$ entre ouverts de $(M,[g])$ se remonte en un diff\'eomorphisme  entre ouverts de $B$ qui pr\'eserve $\omega$.  
Ce diff\'eomorphisme pr\'eserve les champs $\omega$-constants, donc si $f \in \Conf(M)$ et si $(b,\xi)$ appartient au domaine de d\'efinition de $\exp$, alors 
\begin{eqnarray}
\label{equ.invariance-exponentielle}
 f (\exp(b,\xi)) = \exp(f(b),\xi)
\end{eqnarray}

Ainsi le flot local $\phi_X^t$ se remonte sur $B$ et le champ  $X$ sur $M$ se remonte en un champ de vecteurs sur $B$, renot\'e $X$, satisfaisant ${\mathcal L}_X \omega =0$.  Comme l'action de $\phi_X^t$ sur  $B$ pr\'eserve un parall\'elisme, si $X$ n'est pas nul, $X$ est sans singularit\'e sur $B$ (voir \cite[Thm I.3.2]{kobayashi}).
On choisit $b_0 \in B$ au-dessus de la singularit\'e $x_0$ de $X$.  Le flot local $\phi_X^t$ agit sur $B$ en pr\'eservant la fibre de $b_0$. Il existe donc pour tout $t \in \BR$ un \'el\'ement $h^t \in P$ tel que $\phi_X^t.b_0. h^{-t}=b_0$.  Il n'est pas difficile de v\'erifier que $\{h^t\}$ constitue un groupe \`a un param\`etre de $P$, que l'on appelle le \emph{flot d'holonomie de $X$ en $x_0$, relativement \`a $b_0$}.  Changer $b_0$ en $b_0.p$ revient simplement \`a conjuguer le flot d'holonomie en $\{ph^t p^{-1}\}$.  Le flot $\{ h^t \}$ d\'efinit un champ conforme sur $\Ein^{p,q}$, avec une singularit\'e en $o$, que l'on appelle {\it champ d'holonomie de $X$ en $x_0$, relativement \`a $b_0$}. On note ce champ d'holonomie $X_h$.

\begin{remarque}
\label{rem.cas-plat}
Dans le cas o\`u $X$ est un champ de vecteurs conforme sur un voisinage $V$ de $o$ dans $\Ein^{p,q}$, alors on a simplement $X_h=X$. C'est une cons\'equence du fait que le fibr\'e normal de Cartan est la pr\'eimage de $V$ dans $(\Ein^{p,q},PO(p+1,q+1),\omega_G)$, et du Th\'eor\`eme de Liouville.  On obtient qu'un champ de vecteurs sur une vari\'et\'e conform\'ement plate est toujours localement conjugu\'e, au voisinage d'une singularit\'e, \`a son champ d'holonomie.  La conjugaison est de plus un diff\'eomorphisme conforme.
\end{remarque}

Soit maintenant $U$ un ouvert de $M$ et $(f_k : U \rightarrow M)$ une suite de plongements conformes.  Soit $x \in U$ et supposons que la suite $f_k(x)$ converge vers  $y$ lorsque $k \rightarrow \infty$.  On dit  qu'une suite  $(h_k)$ de $P$ est une {\it suite d'holonomie de $(f_k)$ en $x$} s'il existe une suite $(b_k)$ dans la fibre de $x$, contenue dans un compact de $B$, et  telle que $f_k(b_k).h_k^{-1}$ soit \'egalement contenue dans un compact de $B$.
 Si $(g_k)$ est une autre suite d'holonomie de $(f_k)$ en $x$, il r\'esulte de la propret\'e de l'action de $P$ sur $B$ que l'on a $h_k = c_k g_kd_k \ \forall k$, o\`u $(c_k)$ et $(d_k)$ sont deux suites born\'ees de $P$. On dit alors que $(h_k)$ et $(g_k)$ sont {\it \'equivalentes dans $P$}. R\'eciproquement, si $(g_k)$ est une suite d'holonomie de $(f_k)$ en $x$, alors toute suite $(h_k)$ \'equivalente \`a $(g_k)$ dans $P$  est encore une suite d'holonomie de $(f_k)$ en $x$.





 \subsection{Stabilit\'e, semi-compl\'etude, et calcul d'holonomie}
 \label{sec.rappel-stabilite}

Nous rappelons ici la notion de stabilit\'e dynamique, qui va jouer un r\^ole fondamental pour la suite.
Soit $U$ un ouvert de $M$ et $(f_k: U \to M)$ une suite de plongements conformes.  On dit que la suite $(f_k)$ est {\it stable en $x \in U$} si pour toute suite $(x_k)$ de $U$ qui converge vers $x$, la suite $f_k(x_k)$ converge vers une limite $x_{\infty}$ ind\'ependante de $(x_k)$.  La suite $(f_k : U \to M)$ est {\it fortement stable en $x\in U$}
 s'il existe un voisinage $V$ de $x$  ainsi qu'un point $x_{\infty} \in M$ de sorte que $f_k(\overline{V}) \to x_{\infty}$.  Une caract\'erisation de la stabilit\'e et de la stabilit\'e forte en termes d'holonomie a \'et\'e donn\'ee en \cite[Lemme 4.3]{frances-degenere}  : {\it la suite $(f_k)$ est stable en $x$ si et seulement s'il existe une suite d'holonomie $(h_k)$ de $(f_k)$ en $x$ qui soit dans $A^+$. Autrement dit, $h_k$ est de la forme $h_k= \mbox{diag}(\lambda_1(k),\ldots,\lambda_n(k)) \in \CO(p,q)$ avec $\lambda_1(k)\geq \ldots \geq \lambda_n(k) \geq 1$; en particulier les suites $1/\lambda_i(k)$ sont born\'ees. La suite $(f_k)$ est fortement stable si et seulement si $1/\lambda_i(k) \to 0$} pour tout $i$.  Des suites d'holonomies comme ci-dessus sont appel\'ees \emph{stables} et \emph{fortement stables}, respectivement.

Lorsque le flot d'holonomie $\{h^t\}$ d'un champ de vecteurs $X$ a la propri\'et\'e de contracter un segment g\'eod\'esique $[oz]$, et d'\^etre stable en $z$,   on obtient  des propri\'et\'es de semi-compl\'etude du flot $\{ \phi_X^t \}$ sur tout un ouvert~:

\begin{proposition}[voir les Prop 6.1, Prop 6.3 de \cite{frances-riemannien}]
\label{prop.stabilite}
Soit $X $ un champ de vecteurs conforme ayant une singularit\'e en   $x_0$,  avec flot d'holonomie $\{ h^t \}$ en $x_0$ relativement \`a $b_0$.  Il existe $R_0 > 0$ tel que si $\alpha : [0,1] \rightarrow M$ est un  segment g\'eod\'esique conforme issu de $x_0$, dont le d\'eveloppement $\beta = \mathcal{D}_{x_0}^{b_0}(\alpha)$ satisfait \`a :

\begin{itemize}
\item $h^t.[\beta] \subset B(o,R_0)$ pour tout $t \geq 0$,
\item $h^t.[\beta] \to o$ lorsque $t \to \infty$,
\item $(h^{t_k})$ est stable en $\beta(1)$ pour tout $t_k \rightarrow \infty$,
\end{itemize}

alors 
il existe un voisinage $V$ de $\alpha(1)$ tel que :
\begin{enumerate}
\item Le flot $\phi^t_X$ est d\'efini sur $V$ pour tout $t \geq 0$.
\item Le flot $\phi^t_X$ est d\'efini sur $[\alpha]$ pour tout $t \geq 0$, et $\lim_{t \to \infty} \phi_X^t.[\alpha] \to x_0$.
\item Il existe $s_k \rightarrow \infty$ telle que $(\phi^{s_k}_X)$ soit stable en $y$ pour tout $y \in V$, et toute suite d'holonomie de $(h^{s_k})$ en $\beta(1)$ soit suite d'holonomie de $(\phi^{s_k}_X)$ en $y$.  De plus, si $(h^{t_k})$ est fortement stable en $\beta(1)$ pour tout $t_k \rightarrow \infty$, alors la suite $(\phi^{s_k}_X)$ est fortement stable en chaque $y \in V.$
\end{enumerate}
\end{proposition}

Dans le cas de cette proposition o\`u toutes les suites $(h^{t_k})$ sont fortement stables, on peut trouver une suite $( \phi_X^{s_k})$ qui contracte un ouvert $V$ sur $x_0$.   Mais pour montrer l'annulation du tenseur de Weyl sur $V$, il est n\'ecessaire, au vu de l'exemple donn\'e en  Section \ref{sec.exemple-contraction}, d'obtenir des informations plus fines, notamment des informations sur les ``vitesses de contraction"  des diff\'erentielles $D_x\phi_X^t$.  C'est le but du prochain r\'esultat, le Th\'eor\`eme \ref{thm.integrabilite-stable}.

Soient $U$ un ouvert connexe de $M$ et $(f_k : U \to M)$ une suite de plongements conformes.  On suppose que $(f_k)$ est stable en un point $x \in U$. La suite $(f_k)$ admet donc une suite d'holonomie $(h_k)$ qui est dans $A^+$.
\'Ecrivons $h_k=\text{diag}(\lambda_1(k), \ldots , \lambda_n(k))$, avec $\lambda_1(k) \geq \ldots \geq \lambda_n(k) \geq 1$. Il existe $s$ entiers naturels non nuls $n_1, \ldots, n_s$ tels que pour tout $l \in \{0, \ldots, s-1 \}$, et tout couple d'indices $n_{l}+1 \leq i \leq j \leq n_{l+1}$, le quotient $\frac{\lambda_i(k)}{\lambda_j(k)}$ soit born\'e dans $[1, \infty)$  (on a adopt\'e la convention $n_0=0$). Quitte \`a remplacer $(h_k)$ par une suite \'equivalente de ${P}$, on peut alors supposer que pour tout $0 \leq l \leq s-1$, et tout couple d'indices $n_{l}+1 \leq i \leq j \leq n_{l+1}$, $\lambda_i(k)=\lambda_j(k)$ pour tout $k \in \BN$.   Si $j \in \{ 1, \ldots , s \}$, on note alors $\mu_j(k)=\frac{1}{\lambda_{n_j}(k)}$. On a $\mu_1(k) \leq \ldots \leq \mu_s(k) \leq 1$, et l'action de $\Ad h_k$ sur $\lien^-$ se fait par une transformation diagonale de valeurs propres $\mu_1(k), \ldots , \mu_s(k)$. Quitte \`a extraire encore, on peut supposer  que chaque suite  $(\mu_j(k))$ admet une limite dans $\BR^+$, et que $\frac{\mu_{j+1}(k)}{\mu_j(k)} \to \infty$, pour tout $j \in \{1,\ldots,s-1 \}$.  Notons que  les suites $\mu_j(k)$ tendent vers $0$, sauf \'eventuellement pour $j=s$.
   Nous allons maintenant expliciter l'information dynamique que rev\^etent les suites $\mu_j(k)$.

On munit la vari\'et\'e pseudo-riemannienne $(M,g)$ d'une m\'etrique riemannienne auxiliaire, 
qui d\'efinit une norme $|| \cdot ||$ sur $TM$ et une distance $d$ sur $M$.  
L'\'enonc\'e qui suit est ind\'ependant de ce choix d'une m\'etrique auxiliaire.  Rappelons que si $(a_k)$ et $(b_k)$ sont deux suites r\'eelles positives, la notation $a_k=\Theta(b_k)$ signifie qu'il existe $C_1,C_2>0$ tels que pour $k$ suffisamment grand, $C_1 a_k \leq b_k \leq C_2 a_k$.  Par $b_k=O(a_k)$, on entend qu'il existe $C>0$ tel que pour $k$ suffisamment grand, on a $C b_k \leq a_k$.

\begin{theoreme}\cite[Th\'eor\`emes 1.1 et 1.4]{frances-degenere}
\label{thm.integrabilite-stable}
Soit $(f_k)$ une suite stable de transformations locales conformes, d\'efinies sur un ouvert $U$.  Quitte \`a remplacer $(f_k)$ par une suite extraite et \`a r\'etr\'ecir $U$, il existe une filtration de $TU$, ${\mathcal F}_{0}=\{ 0 \} \subsetneq {\mathcal F}_{1} \subsetneq \ldots  \subsetneq {\mathcal F}_{s-1}$, qui s'int\`egre en $s$ feuilletages de $U$, $F_{0} \subsetneq \ldots \subsetneq F_{s-1}$, satisfaisant les propri\'et\'es suivantes~:

\begin{enumerate}
\item{ 
\begin{enumerate}
\item{Un vecteur non nul $u \in T_xU$ appartient \`a  $T_xU \setminus {\mathcal F}_{s-1}(x)$, si et seulement si   pour toute suite $(u_k)$ de $T_xU$ qui converge vers $u$, 
$$||D_xf_k(u_k)|| = \Theta(  \mu_s(k)).$$}
\item{   Un vecteur non nul $u \in T_xU$ appartient \`a ${\mathcal F}_{j}(x) \setminus {\mathcal F}_{j-1}(x)$, $j=1, \ldots s-1$, si et seulement si les deux conditions ci-dessous sont satisfaites~:
\begin{enumerate}
\item{Pour toute suite $(u_k)$ de $T_xU$ qui converge vers $u$,  
$$ \mu_{j}(k)=O(||D_xf_k(u_k)||).$$} 
\item{   Il existe une suite $(u_k)$ de $T_xU$ qui converge vers $u$ telle  que $ ||D_xf_k(u_k)||  = \Theta(\mu_j(k))$.}
\end{enumerate}}
\end{enumerate} }
\item{Chaque $x \in U$ admet un voisinage $U_x$ tel que la feuille locale $F_j^{loc}(x)$ de $x$ dans $U_x$ soit caract\'eris\'ee par :

\begin{enumerate}
\item{Un point $y$  appartient \`a $U_x \setminus  F_{s-1}^{loc}(x)$ si et seulement si pour toute suite $(y_k)$ de $U_x$ qui converge vers $y$,  
$$d(f_k(x),f_k(y_k)) = \Theta(\mu_s(k)).$$}
\item{   Un point $y$ appartient \`a  $ F_{j}^{loc}(x) \setminus F_{j-1}^{loc}(x)$, $j=1, \ldots s-1$, si et seulement si les deux conditions ci-dessous sont satisfaites~:
\begin{enumerate} 
\item{Pour toute suite $(y_k)$ de $U_x$ qui converge vers $y$, 
$$  \mu_{j}(k)=O(d(f_k(x),f_k(y_k))).$$} 
 \item{Il existe une suite $(y_k)$ de $U_x$ qui converge vers $y$ telle que 
 $$d(f_k(x),f_k(y_k)) = \Theta(\mu_j(k)).$$}
 \end{enumerate}}
 \end{enumerate} }
\end{enumerate}

\end{theoreme}

\begin{remarque}
\label{rem.orthogonalite}
Rappelons bri\`evement comment les espaces ${\mathcal F_j(x)}$ sont construits. On d\'ecompose $\lieg/\liep$ en somme directe $ \lieg/\liep= \oplus_{j=1}^s \mathcal{N}_j$
o\`u $(\Ad h_k)_{|{\mathcal N}_j}=\mu_j(k)Id_{{\mathcal N}_j}$.   On d\'efinit ${\mathcal E}_j:= \oplus_{i \leq j} {\mathcal N}_i$ pour $j=\{1, \ldots ,s \}$.  Alors, il existe  $b \in B$  dans la fibre de $x$ pour lequel ${\mathcal F}_j(x)=\iota_b({\mathcal E}_j)$. Ainsi, si $s \geq 3$, on a pour tout $x \in U$, ${\mathcal F}_{s-1}(x)={\mathcal F}_{1}^{\bot}(x)$ (l'orthogonal est pris pour n'importe quelle m\'etrique de la classe conforme).
\end{remarque}

\begin{remarque}
\label{rem.unicite}
Dans \cite{frances-degenere}, Section 8, on a montr\'e que les suites et les feuilletages du Th\'eor\`eme \ref{thm.integrabilite-stable} sont uniques~: s'il existe $r$ suites $(\nu_1(k)), \ldots ,$ $(\nu_r(k))$, avec $\nu_j(k)=o(\nu_{j+1}(k))$, et des sous-vari\'et\'es $H_0^{loc}(x):=\{x \} \subsetneq H_1^{loc}(x) \subsetneq \ldots \subsetneq H_{r-1}^{loc}(x) \subsetneq U_x$ qui satisfont aux conclusions $2(a)$ et $(b)$ du Th\'eor\`eme \ref{thm.integrabilite-stable}, alors n\'eces\-sairement, $r=s$; $\mu_j(k) =\Theta(\nu_j(k) )$; et $H_{j-1}^{loc}(x)=F_{j-1}^{loc}(x)$, pour $j=\{1, \ldots, s\}$.
\end{remarque}

L'int\'er\^et pratique du th\'eor\`eme est qu'il donne des estimations m\'etriques sur la vitesse \`a laquelle les orbites de $(f_k)$ se rapprochent, \`a partir de la connaissance des suites $\mu_j(k)$. R\'eciproquement, en utilisant la propri\'et\'e d'unicit\'e mentionn\'ee ci-dessus, on peut retrouver les suites $\mu_j(k)$ en calculant \`a quelle vitesse les orbites de $(f_k)$ se rapprochent, ce qui nous permet de d\'eterminer rapidement des suites d'holonomies. Nous allons illustrer ceci dans les Propositions \ref{prop.piege} et \ref{prop.piege-lumiere}, o\`u nous allons calculer l'holonomie de $\{ h^t \}$ en des points stables pr\`es de $o$.

Pour donner au lecteur une id\'ee simple de la mani\`ere dont on peut combiner la Proposition \ref{prop.stabilite} et le Th\'eor\`eme \ref{thm.integrabilite-stable}, nous finissons cette section en montrant~:
\begin{proposition}
\label{prop.annulation}
Soit $(M,g)$ une vari\'et\'e pseudo-riemannienne de dimension $\geq 3$, et $X \in \vecconf$ admettant une singularit\'e $x_0$, avec flot d'holonomie $\{h^t\}$ satisfaisante \`a toutes les hypoth\`eses de la Proposition \ref{prop.stabilite}.
\begin{enumerate}
\item{ Si pour toute suite  $t_k \to \infty$, l'holonomie de $(h^{t_k})$ en $\beta(1)$ admet une sous-suite  \'equivalente dans $P$ \`a $(\diag(\lambda_k,\ldots,\lambda_k)) \in (A^+)^{\BN}$, o\`u $\lambda_k \to \infty$, alors il existe un ouvert non vide contenant $\alpha(1)$ qui est conform\'ement plat. }
\item{Si $(M,g)$ est lorentzienne et analytique, et si pour toute suite $t_k \to \infty$,  $(h^{t_k})$ est fortement stable en $\beta(1)$, alors $(M,g)$ est conform\'ement plate.}
\end{enumerate}
\end{proposition}

\begin{proof}
Commen\c{c}ons par la preuve du premier point.  La Proposition \ref{prop.stabilite} donne l'existence d'un ouvert $V$ contenant $\alpha(1)$, et d'une suite $s_k \to \infty$ de sorte que $\phi_X^{s_k}$ est d\'efini sur $V$ pour tout $k$, et $(\phi_X^{s_k})$ admet pour holonomie $(h_k) = (\diag(\lambda_k,\ldots,\lambda_k))$ en chaque point de $V$, avec de plus $\lim_{k \to \infty} \lambda_k=\infty$. En particulier  $(\phi_X^{s_k})$ est fortement stable en chaque point de $V$. On va alors montrer que le tenseur de Weyl  associ\'e \`a $(M,g)$ s'annule sur $V$.  Pour cela, nous utilisons, quitte \`a extraire une nouvelle sous-suite, les conclusions du Th\'eor\`eme \ref{thm.integrabilite-stable}.  Ici, l'entier  $s$ vaut $1$ car  toutes les suites constituant la diagonale de $(h_k)$ sont \'equivalentes, en tant que suite de $\BR$. On a ainsi  $\mu_1(k) = 1/\lambda_k$, qui tend vers $0$ lorsque $k \to \infty$.  La stratification donn\'ee par le th\'eor\`eme est triviale dans ce cas.  Si l'on s'est fix\'e une m\'etrique auxiliaire $\lambda$ sur $M$, et si $|| \cdot ||$ d\'esigne la norme 
 que $\lambda$ d\'efinit sur $TM$, le point $1(a)$ du Th\'eor\`eme \ref{thm.integrabilite-stable} affirme que pour tout $y\in V$, et $u \in T_yM$ non nul, on doit avoir :
$$ ||D_y\phi_X^{s_k}(u)|| = \Theta( \mu_1(k)).$$
 Comme $(\phi_X^{s_k})$ est fortement stable en chaque point de $V$, on peut supposer, quitte \`a restreindre $V$, que l'ensemble $\bigcup_{k \geq 0} \phi_X^{s_k}(V)$ est relativement compact dans $M$.  Par cons\'equent, si $W$ d\'esigne le tenseur de Weyl, on doit avoir sur $\bigcup_{k \geq 0} \phi_X^{s_k}(V)$~:
$$ ||{W}(u,v,w)|| \leq C ||u|| \cdot ||v|| \cdot ||w||,$$
pour une certaine constante $C>0$.  Puis, on utilise l'invariance conforme du tenseur de Weyl; pour $y \in V$, $u,v,w \in T_yM$, et $k \geq 0$ :
$$ D_y\phi_X^{s_k}({W}_y(u,v,w))={W}_{\phi_X^{s_k}.y}(D_y\phi_X^{s_k}(u),D_y\phi_X^{s_k}(v),D_y\phi_X^{s_k}(w)).$$
Si ${W}_y(u,v,w)$ n'est pas nul, le point $1(a)$ du Th\'eor\`eme \ref{thm.integrabilite-stable} donne que $ ||D_y\phi_X^{s_k}({W}_y(u,v,w))||=\Theta(\mu_1(k))$,  tandis que 
$$||{W}_{\phi_X^{s_k}.y}(D_y\phi_X^{s_k}(u),D_y\phi_X^{s_k}(v),D_y\phi_X^{s_k}(w))|| =O(\mu_1(k)^3).$$
On aboutit \`a une contradiction.

Le m\^eme raisonnement assure que le tenseur de Cotton est nul si $\mbox{dim } M = 3$.

Nous passons maintenant au second point de la proposition. La vari\'et\'e $(M,g)$ est suppos\'ee lorentzienne et analytique.  La Proposition \ref{prop.stabilite} donne l'existence d'un ouvert $V$ contenant $\alpha(1)$, et d'une suite $s_k \to \infty$ de sorte que $\phi_X^{s_k}$ est d\'efini sur $V$ pour tout $k$, et $(\phi_X^{s_k})$ admet une holonomie fortement stable $(h_k)$ en chaque point de $V$.  On applique le Th\'eor\`eme \ref{thm.integrabilite-stable}.  Les points $1, (a)$ et $1, (b)$ assurent que pour tout $y$ dans $V$, il existe un voisinage ${\mathcal V}$ de $T_yM$ tel que $D_y\phi_X^{s_k}(\mathcal V) \to 0$.  Autrement dit $(\phi_X^{t_k})$ est fortement stable au sens de \cite{frances-causal} et les Propositions 4 et 5 de \cite{frances-causal} assurent que $V$ est conform\'ement plat. Par analyticit\'e, il en va de m\^eme pour $(M,g)$.  \end{proof}

\section{Preuve du Th\'eor\`eme \ref{thm.translations}}
\label{sec.preuve-translations}
On consid\`ere toujours une structure conforme pseudo-riemannienne $(M,[g])$, de signature $(p,q)$, $p+q \geq 3$.  On suppose que $X$ est un \'el\'ement non trivial de $\vecconf$, admettant une singularit\'e en $x_0 \in M$.  On appelle $(M,B, \omega)$ le fibr\'e normal de Cartan, et $\{h^t\}$ le flot d'holonomie associ\'e \`a $X$ en $x_0$, relativement \`a $b_0 \in B$ au-dessus de $x_0$.  Sous les hypoth\`eses du Th\'eor\`eme \ref{thm.translations} la partie lin\'eaire de la transformation affine $h^t$ est semi-simple sur ${\bf C}$.  On montre alors ais\'ement qu'\`a conjugaison pr\`es  par une translation,  $\{h^t \}$ s'\'ecrit comme le produit commutatif de sa partie lin\'eaire, qui est un groupe \`a un param\`etre relativement compact $\{\kappa^t\}$, et d'un flot de translations $\{\tau^t\}$.

Commen\c{c}ons par \'etudier le cas o\`u $\{\tau^t\}$ est trivial.  Dans ce cas, $h^t={\kappa^t}$ pour tout $t \in \BR$.  Par \cite[Cor 6.2]{frances-riemannien}, il existe un voisinage $U$ de $x_0$ sur lequel $X$ est complet, et engendre un flot relativement compact $\{ \phi_X^t \}$ de transformations conformes de $U$.  On est alors dans le cas $(1)$ du Th\'eor\`eme \ref{thm.translations}.  En moyennant une m\'etrique de la classe conforme par $\{  \phi_X^t\}$, on construit $\tilde g$ sur $U$ pour laquelle $X$ est un champ de Killing. 
Il est alors clair que $X$ est lin\'earisable en $x_0$.  

Nous supposerons donc dor\'enavant que $\{\tau^t\}$ n'est pas trivial.  Deux cas vont devoir \^etre consid\'er\'es s\'epar\'ement : celui o\`u $\{ \tau^t \}$ est un flot de translations de type espace ou temps, et celui o\`u $\{ \tau^t \}$ est un flot de translations de type lumi\`ere.


\subsection{Cas o\`u $\{\tau^t \}$ est un flot de translations  espace ou temps}
\label{sec.cas-spacetemps}

On va dans cette section d\'emontrer le th\'eor\`eme suivant, qui pr\'ecise le Th\'eor\`eme \ref{thm.translations} dans ce cas.  Les notations utilis\'ees sont celles de la Section \ref{sec.convexe}.
\begin{theoreme}
\label{thm.translation-partiel}
Soit $(M,g)$ une vari\'et\'e pseudo-riemannienne de signature $(p,q)$, $p+q \geq 3$.  Soit $X \in \vecconf$, s'annulant en $x_0 \in M$.  On suppose que le flot d'holonomie $\{h^t\}$ de $X$ en $x_0$ est le produit commutatif $\{\kappa^t \cdot \tau^t\}$, avec $\{ \kappa^t \}$ relativement compact et $\{ \tau^t \}$ un groupe de translations de type espace ou de type temps.  Alors il existe un voisinage convexe $U$ de $x_0$, et deux ensembles $U^> $ et $U^< $ dont la r\'eunion est  $U^+$ si $\{ \tau^t \}$ est de type espace, $U^-$ si $\{ \tau^t \}$ est de type temps,   avec  les   propri\'et\'es suivantes~:
\begin{enumerate}
\item{Pour tout $x \in U^>$, le flot $\{ \phi_X^t \}$ est d\'efini pour tout $t \geq 0$, et de plus $\lim_{t \to \infty} \phi_X^t.x=x_0$.}
\item{Pour tout $x \in U^<$, le flot $\{ \phi_X^t \}$ est d\'efini pour tout   $t \leq 0$, et de plus $\lim_{t \to - \infty} \phi_X^t.x=x_0$.}
 \item{L'ouvert $(U^+,[g])$ ou $(U^-,[g])$, suivant que $\{ \tau^t \}$ est de type espace ou temps, est conform\'ement plat.}
\end{enumerate}
\end{theoreme}

Pour v\'erifier que les conclusions du Th\'eor\`eme \ref{thm.translation-partiel} impliquent le cas $(2)$ du Th\'eor\`eme \ref{thm.translations}, il reste \`a montrer que $X$ est essentiel sur tout voisinage de $x_0$.  Mais le flot d'holonomie ne fixe aucun point de l'espace de Minkowski $j(\BR^{p,q})$, donc c'est une cons\'equence de \cite[Prop 4.8]{frances-riemannien} (voir aussi \cite[Thm 2.1]{capocci1}). 

Sur $\BR^{p,q}$, on note 
$$<x,y >:=x_1y_n+ x_ny_1+ \cdots + x_py_{q+1}+ x_{q+1}y_p + x_{p+1}y_{p+1}+\cdots +x_qy_q$$
la forme bilin\'eaire associ\'ee \`a $Q^{p,q}$, et pour simplifier l'\'ecriture, on notera $q(x):=Q^{p,q}(x)$ dans les calculs qui vont suivre.  Nous utiliserons \'egalement les notations introduites dans la Section \ref{sec.convexe}.

On pose :
$$e_1^{\prime}:=\frac{1}{\sqrt 2}(e_1-e_n), \ldots, e_p^{\prime}:=\frac{1}{\sqrt 2}(e_p-e_{q+1})$$
$$e_{p+1}^{\prime}:=e_{p+1}, \ldots , e_q^{\prime}:=e_q$$
$$ e_{q+1}^{\prime}:=\frac{1}{\sqrt 2}(e_p+e_{q+1}), \ldots , e_{n}^{\prime}:=\frac{1}{\sqrt 2}(e_1+e_{n}).$$

Quitte \` a conjuguer $\{ h^t \}$ dans $P$, on peut supposer que $\tau^t$ est la translation de vecteur $tv$ avec $v=e_1^{\prime}$ ou $v=e_n^{\prime}.$  
Le groupe $\{ \kappa^t \}$ est contenu dans un compact maximal de $\mbox{Stab}(v) < \CO(p,q)$.  Comme $v$ est de norme non nulle, ce stabilisateur est semi-simple, et un compact maximal est un produit $\OO(p-1) \times \OO(q)$ ou $\OO(p) \times \OO(q-1)$. Donc si l'on conjugue encore par un \'el\'ement de $P$ qui laisse $v$ invariant, on peut supposer
$\{{\kappa}^t\}$ de la forme suivante~:
\begin{itemize}
\item{un flot de $\OO(p,q)$ laissant $\text{Vect}(e_2^{\prime},\ldots,e_p^{\prime})$ et $\text{Vect}(e_{p+1}^{\prime},\ldots,e_n^{\prime})$ stables, dans le cas o\`u $v=e_1^{\prime}$.}
\item{un flot de $\OO(p,q)$ laissant $\text{Vect}(e_1^{\prime},\ldots,e_p^{\prime})$ et $\text{Vect}(e_{p+1}^{\prime},\ldots,e_{n-1}^{\prime})$ stables, dans le cas o\`u $v=e_n^{\prime}$.}
\end{itemize}
En particulier, $\kappa^t$ vu comme transformation lin\'eaire de $\BR^{p,q}$ laisse  la norme $ || \cdot ||$ invariante (voir \ref{sec.normes} pour les notations).  Nous ferons donc dor\'enavant ces hypoth\`eses sur $\{ \tau^t\}$ et $\{ \kappa^t \}$, ce qui revient, rappelons-le, \`a remplacer le point $b_0$ par un $b_0.p$ ad\'equat dans la m\^eme  fibre.  On notera $\epsilon=1$ si $v=e_n^{\prime}$ et $\epsilon=-1$ si $v=e_1^{\prime}$. Si $U$ est un ouvert convexe contenant $o$, alors  $U^{\epsilon}$ d\'esigne $U^+$ si $\epsilon=1$, et $U^-$ si $\epsilon = -1$.

Pour tout $R >0$, nous posons :
$$ V_R^>:=B^{\epsilon}(o,R) \cap \{ z=j^o(x) \in \BR_o^{p,q} \ | \ \epsilon b(v,x) \geq 0 \},$$
$$ V_R^<:=B^{\epsilon}(o,R) \cap \{ z=j^o(x) \in \BR_o^{p,q} \ | \ \epsilon b(v,x) \leq 0 \}.$$
Notons que $B^{\epsilon}(o,R)=V_R^> \cup V_R^<$.

Si  $z=\pi_G(e^{\xi}) \in B(o,R)$ pour $\xi \in {\mathcal B}(0,R)$, nous appelons $[oz]$ le segment g\'eod\'esique d\'efini par :
$$ [oz]:=\{ \pi_G(e^{u \xi}) \ | \ u \in [0,1]  \}.$$

\begin{proposition}
\label{prop.piege}
 Soit $R>0$. Alors :
\begin{enumerate}
\item{Pour tout $z \in {V_R^>}$ et tout $t \geq 0$, ou pour tout $z \in {V_R^<}$ et tout $t \leq 0$, on a $h^t. [oz] \subset B(o,R)$.  De plus, dans le premier cas, $\lim_{t \to \infty} h^t.[oz]=o$, et $\lim_{t \rightarrow - \infty} h^t.[oz] = o$ dans le second cas.}
\item{Soit $z \in V_R^> \cup V_R^<$ et soit $(t_k)$ une suite de nombres positifs si $z \in V_R^>$, et n\'egatifs si $z \in V_R^<$. Supposons que $|t_k| \rightarrow \infty$.  Alors la suite $( h^{t_k} )$ est fortement stable en $z$ et quitte \`a remplacer $(t_k)$ par une suite extraite, l'holonomie de $(h^{t_k})$ en $z$ est $(h_k) = (\text{diag}(t_k^2, \ldots , t_k^2))$.}
\end{enumerate}
\end{proposition}

\begin{proof}
Rappelons que la classe d'holonomie en un point est ``invariante par perturbation compacte."  En particulier $(h^{t_k})=(\kappa^{t_k}.\tau^{t_k})$ est fortement stable en $z$ si et seulement si $(\tau^{t_k})$ l'est, et la classe d'\'equivalence d'holonomies de $(h^{t_k})$ en $z$ est la m\^eme que celle de $(\tau^{t_k})$.  Par ailleurs,  le flot $\{ \kappa^t \}$ laisse invariante toute boule $B(o,R)$, et $\kappa^t.[oz]=[oz_t]$ avec $z_t=\kappa^t.z$. Il est donc suffisant de montrer la proposition pour ${h^t}={\tau^t}$.

Prenons $z:=j^o(x) \in B^\epsilon(o,R)$.  On a alors $||x|| < R$ et $\epsilon q(x) > 0$. Par ailleurs  le segment $[oz]$ est simplement :
$$ [oz] = \{ j^o(ux) \ | \ u \in [0,1]  \}.$$
Dans la suite des calculs, nous utilisons le fait que comme $v=e_1^{\prime}$ ou $e_n^{\prime}$, on a $||v||=1$
 et $b(y,v)=\epsilon <y,v>$ pour tout $y \in \BR^{p,q}$.

 Le groupe \`a un param\`etre de $\OO(p+1,q+1)$ correspondant \`a $\{ \tau^t \}$ s'\'ecrit~:
$$ 
\left(
\begin{array}{ccc}
1 & tv^* & - \epsilon t^2/2 \\
  & I_{p+q} & -tv   \\
  &         & 1
\end{array}
\right)
$$

o\`u $v^* = v^t.J_{p,q}$ dans la base $e_1, \ldots, e_n$.  
 En cordonn\'ees projectives on peut \'ecrire
$$j^o(ux) = [ 1 : ux_1 : \cdots : ux_n : -\frac{u^2q(x)}{2}].$$

L'image $\tau^t.j^o(ux)$  est alors~:
$$[ 1 + tu \langle x,v \rangle + \frac{\epsilon t^2u^2q(x)}{4} : ux_1 + \frac{t u^2q(x)}{2} v_1 : \cdots : ux_n + \frac{t u^2q(x)}{2} v_n : - \frac{u^2q(x)}{2}].$$

Dans la carte $j^o$ ce point est~:
\begin{equation}
\label{eq.formule1}
x(u,t) = \frac{1}{1 + tu\langle x,v \rangle + \frac{\epsilon t^2 u^2 q(x)}{4}} \cdot \left(ux + \frac{t u^2q(x)}{2} v \right).
\end{equation}

On suppose \`a pr\'esent que $z \in V_R^>$ avec $t \geq 0$ ou $V_R^<$ avec $t \leq 0$.  Cela implique $\epsilon t b(x,v) \geq 0$.  Nous voulons montrer que sous ces conditions $\tau^t.[oz] \subset B(o,R)$, et il suffit pour cela de prouver $||x(u,t)||^2-||ux||^2 \leq 0$.  On calcule donc :
\begin{eqnarray}
\label{equ.majoration}
||x(u,t)||^2 & = & \frac{u^2 ||x||^2 + t u^3q(x) b(x,v) + \frac{t^2 u^4q^2(x)}{4}}{(1 + tu \langle x, v \rangle + \frac{\epsilon t^2u^2 q(x)}{4})^2} \\
& = &  \frac{u^2(||x||^2 + q(x)c(u,t))}{(1 + \epsilon c(u,t))^2}
\end{eqnarray}

avec $c(u,t) = ut b(x,v) + \frac{t^2u^2 q(x)}{4}$. Observons que les conditions $\epsilon q(x)>0$ et $\epsilon t b(x,v) \geq 0$ entra\^{\i}nent que $\epsilon c(u,t) \geq 0$. 

On v\'erifie \`a pr\'esent que :
$$ ||x(u,t) ||^2 - ||u x||^2=\frac{\epsilon c(u,t)(  u^2 \epsilon q(x) -2 u^2 ||x||^2)-u^2c^2(u,t)||x||^2}{(1+\epsilon c(u,t))^2}.$$
On obtient bien que cette quantit\'e est n\'egative puisque $\epsilon c(u,t) \geq 0$, et que par ailleurs $2||x||^2 \geq |q(x)| = \epsilon q(x)$.

On veut maintenant montrer que $\lim_{|t| \to \infty}h^t.[oz]=o$ sous les hypoth\`eses du point (1).  L\`a encore, il suffit de consid\'erer le cas ${h^t}={\tau^t}$.  On a donc toujours la condition $\epsilon t b(x,v) \geq 0$.  
Fixons nous $\delta > 0$.  On cherche $T_{\delta}>0$ tel que $| t| \geq T_{\delta}$ implique :
$$ \sup_{u \in [0,1]}|| x(u,t)|| \leq 2\delta. $$ 
Remarquons que puisque $ ||x(u,t) ||^2 - ||u x||^2\leq 0$, alors 
$$\sup_{u \in [0,\frac{\delta}{R}]}|| x(u,t)|| \leq  \delta,$$
et ce pour tout $t \geq 0$ si $z \in V_R^>$ et pour tout $t \leq 0$ si $z \in V_R^<$.
Maintenant, si  $u \in [\frac{\delta}{R},1]$, on tire de (\ref{equ.majoration}) que :
$$ ||x(u,t)||^2 \leq 16\frac{|| x||^2+|t q(x)|\cdot ||x||+\frac{t^2q^2(x)}{4}}{u^4t^4q^2(x)} \leq 16R^4\frac{(||x||+\frac{|t q(x)|}{2})^2}{t^4\delta^4q^2(x)}.$$
Le majorant est ind\'ependant de $u$ et tend vers $0$ lorsque $ | t | \to \infty$.
Ainsi, il existe $T_{\delta}>0$ tel que pour $t \geq T_{\delta}$ :
$$ \sup_{u \in [\frac{\delta}{R},1]}|| x(u,t)|| \leq \delta,$$
ce qui ach\`eve la preuve du premier point de la proposition.

Nous montrons \`a pr\'esent le second point pour $z \in V_R^>$ et $t_k \to \infty$, le cas $z \in V_R^<$ et $t_k \to -\infty$ \'etant similaire. 
On va prendre une suite $(z_k)$ qui tend vers un point proche de $z$, puis on va calculer, dans la carte \`a l'infini,  la vitesse \`a laquelle $\tau^{t_k}.z$ et $\tau^{t_k}.z_k$ se rapprochent.  Le Th\'eor\`eme \ref{thm.integrabilite-stable} nous permettra alors de d\'eterminer l'holonomie en $z$.

On \'ecrit $z=j^o(x)$.   Soit $x_k=x+w_k$, avec $w_k \to w_{\infty}$, et $z_k=j^o(x_k)$.  Dans les calculs qui suivent, nous notons $\tau^{t_k}.z=j^o(y_k)$ et $\tau^{t_k}.z_k=j^o(y_k^{\prime})$. Nous notons \'egalement 
\begin{eqnarray*}
a = \frac{\epsilon q(x)}{4} & \qquad  & a_k=\frac{\epsilon q(x_k)}{4} \\ 
b=<x,v> & \qquad & b_k=<x_k,v>.
\end{eqnarray*}

De l'expression (\ref{eq.formule1}) pour $u=1$, on tire :
\begin{equation}
\label{eq.formule2}
y_k^{\prime}=\frac{1}{1+t_kb_k+a_kt_k^2} \cdot (x_k+2a_k\epsilon t_k v).
\end{equation}
Si $w_{\infty}$ est dans un petit voisinage de $0$, alors du fait que $q(x) \not =0$, la suite $(a_k)$ va admettre une limite non nulle et donc $y_k^{\prime} \to 0$. On obtient $\lim_{k \to \infty} \tau^{t_k}.z_k=o$, ce qui prouve que $(\tau^{t_k})$ est fortement stable en $z$.  

On suppose d\'esormais que $w_{\infty} \not =0$ est dans un petit voisinage de $0$.
Pour des suites $(a_k)$ et $(b_k)$, on entend par $b_k = o(a_k)$ que $\lim_{k \rightarrow \infty} b_k/a_k = 0$.
De l'expression (\ref{eq.formule2}), on tire :
$$ y_k=\frac{2 \epsilon}{t_k}v + \frac{1}{t_k^2} \left( \frac{x}{a}-\frac{2 \epsilon b}{a}v \right) + o\left(\frac{1}{t_k^2}\right)$$
et 
$$ y_k^{\prime}=\frac{2 \epsilon}{t_k}v + \frac{1}{t_k^2} \left( \frac{x_k}{a_k}-\frac{2 \epsilon b_k}{a_k}v \right) + o\left(\frac{1}{t_k^2}\right).$$

On en conclut que quelle que soit la limite $w_{\infty} \not =0$, $||y_k - y_k^{\prime}|| \sim \frac{C}{t_k^2}$, o\`u la  constante $C$ vaut~:
$$ C = || \frac{x}{a}- \frac{x_\infty}{a_\infty} + \left( \frac{2 \epsilon b_\infty}{a_\infty} - \frac{2 \epsilon b}{a} \right) \cdot v ||.$$
Cette constante ne peut pas \^etre nulle.  S'il en \'etait ainsi, alors on aurait
\begin{eqnarray}
\label{eqn.cnotzero-definite}
\frac{x}{a} - \frac{x_\infty}{a_\infty} = \left( \frac{2 \epsilon b}{a}  - \frac{2 \epsilon b_\infty}{a_\infty}\right) \cdot v.
\end{eqnarray}

Si on applique $\langle \cdot , v \rangle$ aux deux c\^ot\'es, on obtient
$$ \frac{b}{a} - \frac{b_\infty}{a_\infty} = \frac{2 b}{a}  - \frac{2 b_\infty}{a_\infty},$$

ce qui donne une contradiction sauf si $b/a = b_\infty/ a_\infty$.  Dans ce cas, les deux c\^ot\'es de (\ref{eqn.cnotzero-definite}) sont nuls.  Alors $w_\infty = c \cdot x$ pour un $c \in \BR$, et 
$$ \frac{4 x}{\epsilon q(x)} = \frac{4(1+c) x}{\epsilon (1+c)^2 q(x)}$$ 

Cela contredit $w_\infty \neq 0$.

Quitte \`a consid\'erer une sous-suite de $(t_k)$, nous pouvons appliquer le Th\'eor\`eme \ref{thm.integrabilite-stable} \`a la suite $(h^{t_k})$. L'estimation $||y_k - y_k^{\prime}|| \sim \frac{C}{t_k^2}$ quel que soit $w_{\infty} \not = 0$, jointe aux  points $2,(a),(b)$ du Th\'eor\`eme \ref{thm.integrabilite-stable}, donnent que la stratification dynamique est triviale~: l'entier $s$ vaut $1$ et $\mu_1(k)=t_k^{-2}$.  La suite $( \mbox{diag}(t_k^2, \ldots, t_k^2))$ de $A^+$ est une suite d'holonomie pour $(\tau^{t_k})$ en $z$.  \end{proof}

Nous pouvons \`a pr\'esent transcrire sur le flot local $\{ \phi_X^t \}$ les informations r\'ecolt\'ees sur $\{ h^t \}$.   Soit $R_0 >0$ comme dans la Proposition \ref{prop.stabilite}.  Soit $R \leq R_0$ tel que l'application $\xi \mapsto \exp(b_0,\xi)$ soit d\'efinie et injective sur $\mathcal{B}(0,R) \subset \mathfrak{n}^-$.  Soient $\mathcal{V}^>_R$ et $\mathcal{V}^<_R \subset B^\epsilon(0,R)$ les images r\'eciproques $(j^o)^{-1}(V_R^>)$ et $(j^o)^{-1}(V_R^<)$, respectivement.  Pour $\xi \in \mathcal{V}_R^>$, si $\alpha(u) = \exp(b_0,u\xi)$, alors la Proposition \ref{prop.piege} nous dit que $\{ h^t \}$ satisfait aux hypoth\`eses de la Proposition \ref{prop.stabilite} avec $\beta = \mathcal{D}_{x_0}^{b_0}(\alpha)$; il en va de m\^eme pour $\xi \in \mathcal{V}_R^<$ et $\{ h^{-t} \}$.  Posons $U^> := \pi(\exp(b_0,\mathcal{V}^>_R))$ et $U^< := \pi(\exp(b_0,\mathcal{V}^<_R))$. La Proposition \ref{prop.stabilite} implique alors que $\phi^t_X$ est d\'efini pour tout temps positif sur $U^>$ et pour tout temps n\'egatif sur $U^<
 $; de plus $\lim_{t \rightarrow  \infty} \phi_X^t.x = x_0$ pour tout $x \in U^>$ et  $\lim_{t \rightarrow  -\infty} \phi_X^t.x = x_0$ pour tout $x \in U^<$, comme annonc\'e dans les points (1) et (2) du Th\'eor\`eme \ref{thm.translation-partiel}.  Par ailleurs, par le point $(2)$ de la Proposition \ref{prop.piege}, et le point $(1)$ de la Proposition \ref{prop.annulation}, il existe un ouvert non vide contenant $\alpha(1)$ qui est conform\'ement plat.  Comme on a choisi $\xi$ quelconque dans $\mathcal{V}_R^> \cup \mathcal{V}_R^<$, on obtient que $U^> \cup U^<$ est conform\'ement plat. Comme $U^> \cup U^<=U^{\epsilon}$, avec $U = \exp(b_0,\mathcal{B}(0,R))$ convexe, on obtient le point $(3)$ du Th\'eor\`eme \ref{thm.translation-partiel}.

\subsection{Cas o\`u $\{ \tau^t \}$ est un flot de translations de type lumi\`ere}
Nous allons, dans cette section, pr\'eciser le Th\'eor\`eme \ref{thm.translations}, en prouvant~:

\begin{theoreme}
\label{thm.translations-partiel-lumiere}
Soit $(M,g)$ une vari\'et\'e pseudo-riemannienne de signature $(p,q)$, $p \geq 1$, $p+q \geq 3$.  Soit $X \in \vecconf$, s'annulant en $x_0 \in M$.  On suppose que le flot d'holonomie $\{h^t\}$ de $X$ en $x_0$ est le produit commutatif $\{ \kappa^t \cdot \tau^t \}$, avec $\{ \kappa^t \}$ relativement compact et $\{ \tau^t \}$ un groupe de translations de type lumi\`ere.  Alors il existe~:
\begin{itemize}
\item un segment g\'eod\'esique de lumi\`ere $\Delta$ passant par $x_0$, et sur lequel $X$ s'annule;
\item deux ouverts $U^>$ et $U^<$, avec $\Delta \subset \overline{U^>} \cap \overline{U^<}$; et
\item deux submersions $\pi^> : U^> \rightarrow \Delta$ et $\pi^< : U^< \rightarrow \Delta$, dont les fibres sont des hypersurfaces d\'eg\'en\'er\'ees,
\end{itemize}
avec le propri\'et\'es suivantes :
\begin{enumerate}
\item Pour tout $x \in U^>$, le flot $\phi^t_X$ est d\'efini pour tout $t \geq 0$, et de plus $\lim_{t \rightarrow \infty} \phi^t_X .x = \pi^>(x).$
\item Pour tout $x \in U^<$, le flot $\phi^t_X$ est d\'efini pour tout $t \leq 0$, et de plus $\lim_{t \rightarrow - \infty} \phi^t_X.x = \pi^<(x)$.
\item L'ouvert $(U,[g])$, o\`u $U = U^> \cup U^<$, est conform\'ement plat.
\end{enumerate}
\end{theoreme}

 L'\'etude de ce cas est proche de celle qui a \'et\'e men\'ee dans \cite{nilpotent}, et ce qui suit est une adaptation au cadre des champs de vecteurs des techniques introduites dans  \cite[Secs 4.1 - 5.1]{nilpotent}. 
 
Quitte \`a conjuguer $\{ h^t \}$ dans $P$, ce qui revient \`a remplacer $b_0$ par un point $b_0.p$, pour $p \in P$ ad\'equat, on peut supposer que $\tau^t$ est la translation de vecteur $tv$ avec $v=e_1$ et que ${ \kappa^t}$ est inclus dans $\mbox{CO}(p,q)$, et  pr\'eserve la norme $|| \cdot ||$ pour tout $t$. 
Nous reprenons les notations des Sections \ref{sec.normes} et \ref{sec.convexe}. 
On v\'erifie que $\Sigma:={\mathcal N}(\lien^-) \cap {\mathcal S}(0,1)$  est une vari\'et\'e lisse diff\'eomorphe \`a ${\bf S}^{p-1} \times {\bf S}^{q-1}$  (si $p=1$, alors $\Sigma$ n'est pas connexe).   Il existe dans $\lien^-$ un vecteur $\xi_1$ caract\'eris\'e par le fait que $\text{Fix}(\Ad \tau^t) \cap \lien^-=\BR. \xi_1$ (\cite[Lem 4.6]{nilpotent}).  La g\'eod\'esique $\Lambda(v):=\pi_G(e^{v \xi_1})$ est de type lumi\`ere, et l'on dit que c'est la g\'eod\'esique singuli\`ere associ\'ee \`a $\{ \tau^t \}$.  On note ${\hat \Lambda}(v):=e^{v\xi_1}$.  Remarquons que puisque $\{ \kappa^t \}$ commute \`a $\{ \tau^t \}$ et $\{\Ad \kappa^t \}$ pr\'eserve $\lien^-$, on doit avoir $(\Ad \kappa^t).\xi_1=\xi_1$ pour tout $t$.  On conclut que $h^t$ fixe $\Lambda$ point par point pour tout $t$.    Pour $\delta>0$, on d\'efinit:
$$ \Lambda_{\delta}:= \{  \Lambda(v) \ | \  v \in (-\delta,\delta)  \}.$$

L'hyperplan $\xi_1^{\bot}$, o\`u l'orthogonal est pris relativement \`a $\lambda^{p,q}$, s\'epare le c\^one ${\mathcal N}(\lien^-) $ en deux composantes connexes $\mathcal S$ et $- \mathcal S $, et l'on va noter  ${\mathcal S}_{\Sigma}=\Sigma \cap {\mathcal S}$.  
Pour tout $r >0$, d\'efinissons :

$$C^>(r) :=\{  \pi_G(e^{ u \xi}) \in \Ein^{p,q} \ | \  u \in (0,r), \ \xi \in {\mathcal S}_{\Sigma} \},$$
et 
$$C^<(r) :=\{  \pi_G(e^{ u \xi}) \in \Ein^{p,q} \ | \  u \in (0,r), \ \xi \in   -{\mathcal S}_{\Sigma} \}.$$

Nous appellerons par la suite  $C(r)$ la r\'eunion $C^>(r)  \cup C^<(r) \cup \{o \}$, et pour $z \in C^>(r)  \cup C^<(r)$, on notera $[oz]$ le segment g\'eod\'esique de lumi\`ere joignant $o$ et $z$.

\begin{proposition}
\label{prop.piege-lumiere}
Pour $r>0$, on a les propri\'et\'es suivantes :
\begin{enumerate}
\item{Pour tout $z \in C^>(r) $  et tout $t \geq 0$, ou pour tout $z \in C^<(r)$ et $t \leq 0$, on a l'inclusion $h^t.[oz] \subset C(r)$ et $\lim_{|t| \to \infty}h^t.[oz]=o$. }
\item{Soit $z \in C(r)^> \cup C(r)^<$ et soit $(t_k)$ une suite de nombres positifs si $z \in C(r)^>$, et n\'egatifs si $z \in C(r)^<$. Supposons que $|t_k| \rightarrow \infty$. Alors la suite $(h^{t_k})$ est  stable en $z$.  Quitte \`a remplacer $(t_k)$ par  une suite extraite,  l'holonomie de $(h^{t_k})$ en $z$ est $h_k=\text{diag}(t_k^2,|t_k|, \ldots , |t_k|, 1)$.}
\end{enumerate}

\end{proposition}

\begin{proof}
Pour les m\^emes raisons que celles qui ont \'et\'e  expliqu\'ees au d\'ebut de la preuve de la Proposition  \ref{prop.piege}, il suffit de faire la preuve lorsque $\{ h^t \}=\{ \tau^t \}$. 
Le groupe \`a un param\`etre de $\OO(p+1,q+1)$ correspondant \`a $\{ \tau^t \}$ s'\'ecrit comme
$$ 
\left(
\begin{array}{ccc}
1 & tw^* & 0 \\
  & I_{p+q} & -tw   \\
  &         & 1
\end{array}
\right)
$$
o\`u $w=e_1$ et $w^*$ est le transpos\'e de $e_n$.
On choisit $z \in C^>(r)  \cup C^<(r)$. Alors $z=j^o(x)$ pour un certain $x=(s_1, \ldots , s_n)$ satisfaisant $q(x)=0$ et $s_n \not =0$. 
 On a $s_n>0$ si $z \in C^>(r)$ et $s_n<0$ si $z \in C^<(r)$.  Le segment $[oz]$ est $\{ \beta(u):=j^o(ux) \ | \ u \in [0,1]  \}.$  Pour tout $t$, on calcule que $\tau^t.\beta(u)=j^o(x(u,t))$ avec :
\begin{equation}
\label{eq.formule4}
x(u,t)=\frac{ux}{1+tus_n}
\end{equation}

On en d\'eduit donc que $\tau^t.\beta(u)=\beta(\frac{u}{1+tus_n})$.  Ainsi $\tau^t.\beta(u) \in C^>(r)$ pour tout $t \geq 0$ si $z \in C^>(r)$, et $\tau^t.\beta(u)\in C^<(r)$ pour tout $t \leq 0$ si $z \in C^<(r)$.  De plus, $\lim_{|t| \to \infty}[\beta]=\beta(0)=o$ sous ces hypoth\`eses.

On va maintenant montrer le second point de la proposition. On suppose $z \in C^>(r) $ et $t_k \to \infty$ (le cas $z \in C^<(r) $ et $t_k \to -\infty$ se traite de mani\`ere identique). Soit $x_k:=(s_1(k), \ldots ,s_n(k))$ une suite qui tend vers $x^{\prime}:=(s_1^{\prime},\ldots,s_n^{\prime})$ assez proche de $x$ pour que $s_n s_n^{\prime} >  0$. Soit  $z_k=j^o(x_k)$.  Dans les calculs qui suivent, nous notons $\tau^{t_k}.z=j^o(y_k)$, $\tau^{t_k}.z_k=j^o(y_k^{\prime})$ et $\sigma_k=\frac{q(x_k)}{2}$. On calcule :
\begin{equation}
\label{eq.formule5}
 y_k^{\prime}= \left( \frac{s_1(k)+t_k\sigma_k}{1+t_ks_n(k)}, \frac{s_2(k)}{1+t_ks_n(k)}, \ldots, \frac{s_n(k)}{1+t_ks_n(k)} \right).
 \end{equation}
Si $x_k \to x$, alors $\sigma_k \to 0$.  On a donc $y_k^{\prime} \to 0$ et $\tau^{t_k}.z_k \to o$ : la suite $(\tau^{t_k})$ est stable en $z$.

On va d\'eterminer son holonomie en $z$, quitte \`a consid\'erer une suite extraite de $(t_k)$.   Pour cela, on va appliquer les conclusions du  Th\'eor\`eme \ref{thm.integrabilite-stable} \`a la suite $(f_k):=(\tau^{t_k})$, qui est stable en $z$.  Il existe un entier $s \geq 1$, un voisinage ouvert $U_z$ de $z$ dans $\Ein^{p,q}$, des suites $\mu_1(k)< \cdots < \mu_s(k)$ et des sous-vari\'et\'es $F_0^{loc}(z) \subsetneq F_1^{loc}(z) \subsetneq \cdots \subsetneq F_{s-1}^{loc}(z) \subsetneq U_z $ qui satisfont aux conclusions $2(a)$ et $(b)$ du th\'eor\`eme.  En gardant les m\^emes notations que ci-dessus, on va supposer que $z^{\prime}:=j^o(x^{\prime}) \in U_z$, et que la distance $d$ intervenant dans les points $2(a)$ et $(b)$ est la pouss\'ee par l'application $j^o$ de la distance induite par la norme $\| \cdot \|$.     Posons~:
$$E_2= j^o(\mathcal{N}^{p,q}) \cap U_z$$
et
$$E_{1}= j^o(\BR x)  \cap U_z.$$

On reprend les formules (\ref{eq.formule4}) et (\ref{eq.formule5}) et l'on obtient~:
$$ y_k=\frac{x}{1+t_ks_n}=\frac{x}{t_ks_n}-\frac{x}{t_k^2s_n^2}+o\left( \frac{1}{t_k^2} \right)$$
et
$$ y_k^{\prime}=\frac{t_k \sigma_k}{1+t_ks_n(k)} \cdot e_1 + \frac{x_k}{t_ks_n(k)}-\frac{x_k}{t_k^2s_n(k)^2}+o \left( \frac{1}{t_k^2} \right).$$

Ces expressions vont nous permettre de comprendre \`a quelle vitesse $y_k$ et $y_k^{\prime}$ se rapprochent.  Dans la suite, on pose $\sigma_{\infty}=\lim_{k \to \infty} \sigma_k$.  
\begin{itemize}
\item{On suppose que $x^{\prime} \not \in \mathcal{N}^{p,q}$.  Alors $\sigma_{\infty} \not = 0$. La suite 
$(||y_k - y_k^{\prime}||)$ tend vers $\frac{\sigma_{\infty}}{s_n^{\prime}}$.  Autrement dit, pour tout $z^{\prime} \in U_z \setminus E_2$, et toute suite $(z_k)$ qui converge vers $z^{\prime}$, il existe $C>0$ tel que $d(\tau^{t_k}.z,\tau^{t_k}.z_k) \sim C$.}

\item{Supposons maintenant que  $x^{\prime} \in \mathcal{N}^{p,q} \setminus \BR x$.  On a $\sigma_{\infty}=0$, et $x^{\prime} \not \in \BR x$.  La composante de la diff\'erence $y_k - y_k^{\prime}$ selon le vecteur $e_j$, pour $j=2, \ldots , n$ est~:
\begin{equation}
\label{equ.difference}
\left( \frac{s_j}{s_n}-\frac{s_j(k)}{s_n(k)} \right) \frac{1}{t_k}+ \left( \frac{s_j(k)}{s_n(k)^2}-\frac{s_j}{s_n^2} \right) \frac{1}{t_k^2}+o\left( \frac{1}{t_k^2} \right).
 \end{equation}
 
Notons que le coefficient de $\frac{1}{t_k}$ est trivial pour $j = n$.  Si l'\'egalit\'e  $s_j=\frac{s_n}{s_n^{\prime}}s_j^{\prime}$ a lieu pour tout $j=2, \ldots ,n$ alors de  $q(x)=q(x^{\prime})$, on tire $s_1=\frac{s_n}{s_n^{\prime}}s_1^{\prime}$---une contradiction avec  $x^{\prime} \not \in \BR x$.  On conclut que si $z^{\prime} \in E_2 \setminus E_1$, alors pour toute suite $(z_k)$ de $U_z$ qui tend vers $z^{\prime}$, $\frac{1}{t_k}=O(d(\tau^{t_k}.z,\tau^{t_k}.z_k))$.
   
   Soit $(x_k)$ une suite qui tend vers $x^{\prime}$.  Alors la suite $(\sigma_k)$ tend vers $0$.  Aussi, quitte \`a remplacer $(x_k)$ par une suite extraite, tout en gardant {\it la m\^eme} suite $(t_k)$, on peut supposer que $\sigma_k=o(\frac{1}{t_k})$.  Dans ce cas, la composante de $y_k - y_k^{\prime}$ selon le vecteur $e_j$, pour $j = 1,\ldots,n$, est de la forme~:
$$\frac{1}{t_k} \left( \frac{s_j}{s_n}-\frac{s_j(k)}{s_n(k)} \right)+o\left( \frac{1}{t_k}\right).$$ 
 L\`a encore, l'hypoth\`ese $x^{\prime} \not \in \BR x$ assure que l'une des limites $\frac{s_j}{s_n}-\frac{s_j'}{s_n^{\prime}}$ n'est pas nulle.  On en conclut que $ \|y_k - y_k^{\prime}\| = \Theta(\frac{1}{t_k})$.  Donc, si $z^{\prime} \in E_2 \setminus E_1$, il existe $(z_k)$ qui converge vers $z^{\prime}$ telle que~:
$$ d(\tau^{t_k}.z,\tau^{t_k}.z_k) = \Theta \left( \frac{1}{t_k} \right).$$}

\item{Soit enfin  $x^{\prime}\in \BR x \setminus \{x \}$.  Autrement dit $x^{\prime}= \lambda x$, avec $\lambda \not = 1$.  Quitte \`a r\'eduire l'ouvert $U_z$, on va supposer que $\lambda \in ]0,2[$.  De l'expression (\ref{equ.difference}), on d\'eduit que la composante selon $e_n$  de la diff\'erence $y_k - y_k^{\prime}$ est de la forme $(\frac{1}{\lambda}-1)\frac{1}{s_n}\frac{1}{t_k^2}+o(\frac{1}{t_k^2})$.  On conclut que si $z^{\prime} \in E_1 \setminus \{z \}$ est proche de $z$, et si $(z_k)$ est une suite qui tend vers $z^{\prime}$, alors~:
$$   \frac{1}{t_k^2}=O(d(\tau^{t_k}.z,\tau^{t_k}.z_k)).$$

Finalement,  consid\'erons une suite quelconque $x_k \rightarrow x'$.  La suite $(\sigma_k)$ tend vers $0$, donc quitte \`a remplacer $(x_k)$ par une sous-suite, sans modifier la suite $(t_k)$, on peut supposer que $\sigma_k$ est n\'egligeable devant $\frac{1}{t_k^2}$.  Dans ce cas, la formule (\ref{equ.difference}) est aussi valide pour $j=1$.  Comme $\frac{s_j}{s_n} - \frac{s_j(k)}{s_n(k)}$ tend vers $0$ pour $j = 1, \ldots, n-1$, on peut \`a nouveau remplacer  $(x_k)$ par une suite extraite, sans modifier $(t_k)$, de sorte que cette suite de diff\'erences soit un $O(\frac{1}{t_k^2})$.  Ainsi, pour $j=1,\ldots,n-1$, toutes les composantes de $y_k-y_k^{\prime}$ selon $e_j$   sont en $O(\frac{1}{t_k^2})$.  Par le paragraphe pr\'ec\'edent, la composante de $y_k-y_k^{\prime}$ selon $e_n$ est \'equivalente \`a $(\frac{1}{\lambda}-1)\frac{1}{s_n}\frac{1}{t_k^2}$.  On en conclut qu'il existe $z_k \rightarrow z^{\prime}$ telle que~:
$$  d(\tau^{t_k}.z,\tau^{t_k}.z_k) = \Theta \left( \frac{1}{t_k^2} \right).$$
} 
\end{itemize}

En conclusion, le Th\'eor\`eme \ref{thm.integrabilite-stable}, et l'unicit\'e mentionn\'ee en Remarque \ref{rem.unicite}, montrent que l'entier $s$ vaut $3$, que $\mu_1(k)=\frac{1}{t_k^2}$, $\mu_2(k)=\frac{1}{t_k}$, et $\mu_3(k)=1$.  Enfin, $F_{1}^{loc}(z)=E_1$ et $F_{2}^{loc}(z)=E_2$.

 Par la mani\`ere dont les suites $(\mu_j(k))$ sont d\'efinies \`a partir de l'holonomie, on conclut que $(h_k)=(\text{diag}(t_k^2, t_k, \ldots , t_k, 1))$ est une suite d'holonomie de $(\tau^{t_k})$ en $z=j(x)$, comme annonc\'e. \end{proof}

Sur la vari\'et\'e $M$, on pose ${\hat \Delta}(v):=\exp(b_0,v\xi_1)$, d\'efinie sur un intervalle ouvert $I$ contenant $0$, et $\Delta(v)=\pi({\hat \Delta}(v)).$  Pour tout $\delta>0$ tel que $(-\delta,\delta) \subset I$, on pose $\Delta_{\delta}:= \{ \Delta(v) \ | \ v \in (-\delta, \delta)  \}$.
Du fait que $(\Ad h^t).\xi_1=\xi_1$ pour tout $t$, on obtient par exactement la m\^eme preuve que \cite[Prop 5.5]{nilpotent}  :

\begin{lemme}
\label{lem.invariance.holonomie}
Pour tout $v \in (-\delta, \delta)$, le flot d'holonomie de $X$ en $\Delta(v)$ relativement \`a $\hat \Delta(v)$ est $\{ h^t \}$.
\end{lemme}

On d\'efinit pour tout $\delta>0$ tel que $(-\delta,\delta) \subset I$, et tout $r>0$ assez petit une application 
\begin{eqnarray*}
\psi_1 & : & ( -\delta, \delta) \times (0,r) \times ( {\mathcal S}_{\Sigma} \cup -{\mathcal S}_{\Sigma}) \to M \\
& & (v,u,\xi) \mapsto \pi(\exp({{\hat \Delta}(v), u \xi})).
\end{eqnarray*}

L'image de $\psi_1$ est not\'ee $\dot{U}(\delta,r)$.  On d\'efinit les sous-ensembles :
$$\dot{U}^>(\delta,r):=\psi_1((-\delta, \delta) \times (0,r) \times {\mathcal S}_{\Sigma} ),$$ 
et 
$$\dot{U}(\delta,r)^<: = \psi_1((-\delta, \delta) \times (0,r) \times -{\mathcal S}_{\Sigma} ).$$ 
 Chacun de ces ensembles  contient $\Delta_{\delta}$ dans son adh\'erence.  On note  ${U}(\delta,r)=\dot{U}(\delta,r) \cup \Delta_{\delta}$,  et les analogues pour ${U}^>(\delta,r)$ 
 et $U^<(\delta,r)$.


On va maintenant trouver $\delta_0$ et $r_0>0$ tels que  $\phi_X^t.x$ est d\'efini pour tout $t \geq 0$ si $x \in {U}^>(\delta_0,r_0)$, et pour tout $t \leq 0$ si $x \in {U}^<(\delta_0,r_0)$.


On choisit $\delta_0>0$ et $r_0>0$ suffisamment petits pour que 
\begin{itemize}

\item $\xi \mapsto \pi(\exp(\hat \Delta(v),\xi))$ soit une injection sur ${\mathcal B}(0,r_0)$ pour tout $v \in (- \delta_0, \delta_0)$
\item $C( r_0) \subset B(o,R)$, o\`u $R >0$ est donn\'e par la Proposition \ref{prop.stabilite}
\end{itemize}

   On \'ecrit $x=\pi(\exp(\hat \Delta(v),u \xi))$, avec $\xi \in {\mathcal S}_{\Sigma}$. On d\'efinit \'egalement~:
   $$\alpha(s):=\pi(\exp(\hat \Delta(v),su \xi))) \text{ et } \beta:={\mathcal D}_{\Delta(v)}^{\hat \Delta(v)}(\alpha).$$ Il s'agit d'un segment g\'eod\'esique de lumi\`ere de  $C^>(r_0) \subset B(o,R)$.  Le flot d'holonomie de $\{ \phi^t_X \}$ en $\Delta(v)$ relativement \`a $\hat \Delta(v)$ est $\{ h^t \}$ par le Lemme \ref{lem.invariance.holonomie}. On peut alors appliquer  les Propositions \ref{prop.piege-lumiere}  et \ref{prop.stabilite}.  On conclut que  $\phi_X^t.\alpha(1)$ est d\'efini pour tout $t \geq 0$ et que $\lim_{t \to \infty} \phi_X^t.x=\alpha(0)=\Delta(v)$.   Le cas $\xi \in -{\mathcal S}_{\Sigma}$ se traite de mani\`ere identique. 
   
  Soit de plus $\delta_0$ assez petit pour que $\Delta$ soit injective sur $(- \delta_0, \delta_0)$.  Dans ce cas, on va montrer que $\psi_1$ est injective sur $( -\delta_0, \delta_0) \times (0,r_0) \times ( {\mathcal S}_{\Sigma} \cup -{\mathcal S}_{\Sigma})$.   Soit $x \in \dot{U}^>(\delta_0,r_0)$ avec $x=\psi_1(v,u,\xi)=\psi_1(v^{\prime},u^{\prime},\xi^{\prime})$. Comme $\lim_{t \to \infty} \phi_X^t.x=\Delta(v)=\Delta(v^{\prime})$, alors n\'ecessairement $v=v^{\prime}$.  Puis comme $\zeta \mapsto \pi(\exp(\hat \Delta(v),\zeta))$ est  injective sur ${\mathcal B}(0,r_0)$, on obtient $u=u^{\prime}$ et $\xi=\xi^{\prime}$.  Par le Th\'eor\`eme d'Invariance du Domaine,  les ensembles $\dot{U}^>(\delta_0,r_0)$ et $\dot{U}^<(\delta_0,r_0)$ sont ouverts.
   
Alors soit $U^> = \dot{U}^>(\delta_0,r_0)$ et $U^< = \dot{U}^<(\delta_0,r_0)$ pour $\delta_0$ et $r_0$ comme ci-dessus.  On d\'efinit $\pi^>$ sur $U^>$ par 
$$ \pi^> : \psi(v,u,\xi) \mapsto \Delta(v)$$
et similairement pour $\pi^<$ sur $U^<$.    Ces ensembles et projections satisfont \`a (1) et \`a (2) du Th\'eor\`eme \ref{thm.translations-partiel-lumiere}, en prenant $\Delta_\delta$ pour le segment lumi\`ere.   Pour (3), il suffit de prouver la proposition suivante. 
   
 
  
 \begin{proposition}
\label{prop.plat-lumiere}
La m\'etrique est conform\'ement plate sur $U^>$ et $U^<$.
\end{proposition}

\begin{proof}
Nous allons montrer que la g\'eom\'etrie est conform\'ement plate sur $U ^> = \dot{U}^>(\delta_0,r_0)$. Le cas de $U^<$ se traite de mani\`ere identique.  

Soit donc $x \in \dot{U}^>(\delta_0,r_0)$, que l'on \'ecrit $x=\psi_1(v,u,\xi)$.  
Soit $\alpha$ et $\beta$ comme dans la preuve pr\'ec\'edente.  Le Lemme \ref{lem.invariance.holonomie} assure que le flot d'holonomie de $X$ en $\Delta(v)$ est $\{ h^t \}$. Par les Propositions \ref{prop.stabilite} et \ref{prop.piege-lumiere}, il existe  un voisinage $U^{\prime} \subset \dot{U}^>(\delta_0,r_0)$ de $x$
 et une suite $s_k \to \infty$ telle que $(\phi_X^{s_k})$ soit stable en tout point $y \in U^{\prime}$ avec  pour holonomie $(h_k):=(\text{diag}(s_k^2,s_k,\ldots,s_k,1))$.  

Quitte \`a remplacer $(s_k)$ par une suite extraite, on peut appliquer le Th\'eor\`eme \ref{thm.integrabilite-stable}. Les suites intervenant dans le th\'eor\`eme sont ici $\mu_1(k)=1/s_k^{2}$, $\mu_2(k)=1/s_k$ et $\mu_3(k)=1$. Il en r\'esulte deux feuilletages $F_1$ et $F_2$. Le feuilletage  $F_1$ consiste en des courbes de type lumi\`ere. Par la Remarque \ref{rem.orthogonalite}, $\mathcal{F}_2(x) = \mathcal{F}_1^\perp(x)$ pour tout $x \in U'$, si bien que $F_2$ consiste en des hypersurfaces de type lumi\`ere. 

D\'esignons par $C_{{\Delta}(v)}$ le c\^one de lumi\`ere local issu de ${\Delta}(v)$.  Les points $2(a), (b)$ du Th\'eor\`eme \ref{thm.integrabilite-stable} nous disent que pour $z \in U^{\prime}$, l'espace ${\mathcal F}_2(z)$ est tangent \`a la vari\'et\'e stable locale de $z$, qui est constitu\'ee des points $y$ proches de $z$ tels que $\lim_{k \to \infty}\phi_X^{s_k}.z=\lim_{k \to \infty}\phi_X^{s_k}.y$. Il ressort du point (1) d\'ej\`a prouv\'e du Th\'eor\`eme \ref{thm.translations-partiel-lumiere} que ${\mathcal F}_2(z)$ est tangent \`a $C_{{\Delta}(v)} \cap U^{
 \prime}$ en $z$. La Remarque \ref{rem.orthogonalite} assure que ${\mathcal F}_2(z)={\mathcal F}_1^{\bot}(z)$.  Les points $1(a)$ et $1(b)$ du Th\'eor\`eme \ref{thm.integrabilite-stable} assurent qu'en chaque $z \in U^{\prime}$, il existe un rep\`ere $(u_1(z), \ldots ,u_n(z))$ de $T_zM$, avec $u_1(z), \ldots, u_{n-1}(z) \in {\mathcal F}_2(z)$ et $u_1(z) \in {\mathcal F}_1(z)$, ainsi qu'une suite de rep\`eres $(u_{1,k}(z), \ldots, u_{n,k}(z))$, o\`u $u_{i,k}(z) \to u_i(z)$, pour  $i=1, \ldots ,n$, satisfaisant :
\begin{equation}
\label{eq1}
 ||D_z \phi_X^{s_k}(u_{1,k}(z))|| =\Theta \left( \frac{1}{s_k^2} \right).
 \end{equation}
 \begin{equation}
 \label{eq2}
||D_z \phi_X^{s_k}(u_{j,k}(z))|| =\Theta \left( \frac{1}{s_k} \right), \ \ j=2, \ldots ,n-1.
\end{equation}
\begin{equation}
\label{eq3}
||D_z \phi_X^{s_k}(u_{n,k}(z))|| =\Theta( 1).
\end{equation}

La preuve est maintenant une petite perturbation de celle de \cite[Prop 6.1]{nilpotent}. 
Nous en r\'eexpliquons bri\`evement l'id\'ee lorsque $\text{dim }M \geq 4$  (le cas de la dimension $3$ est similaire en rempla\c{c}ant le tenseur de Weyl par le tenseur de Cotton). Soit $W$ le tenseur de Weyl sur $M$. Pour  $U^{\prime}$ suffisamment petit~:
\begin{equation}
\label{eq0}
 ||{W}(u,v,w)|| \leq C ||u|| \cdot ||v|| \cdot ||w||,
 \end{equation}
pour une certaine constante $C>0$. Nous voulons montrer que $W_z=0$ d\`es que $z \in U^{\prime}$.  Pour tout $z \in U^{\prime}$, et tout triplet $(i,j,l) \in \{1, \ldots ,n \}^3$, on \'ecrit~:

\begin{eqnarray*}
 & & D_z\phi_X^{s_k}(W(u_{i,k}(z),u_{j,k}(z),u_{l,k}(z)))= \\
& &\qquad  W_{\phi_X^{s_k}.z}(D_z\phi_X^{s_k}(u_{i,k}(z)), D_z\phi_X^{s_k}(u_{j,k}(z)), D_z\phi_X^{s_k}(u_{l,k}(z))).
 \end{eqnarray*}

De (\ref{eq0}), des estimations (\ref{eq1}), (\ref{eq2}), (\ref{eq3}) et de la relation pr\'ec\'edente, on tire que si $(i,j,l) \not = (n,n,n)$, alors 
$D_z\phi_X^{s_k}(W(u_{i,k}(z),u_{j,k}(z),u_{l,k}(z)))=o(1)$. Des  points $1(a)$ et $1(b)$ du Th\'eor\`eme \ref{thm.integrabilite-stable}, on conclut que  $\text{Im }W$ est inclus dans ${\mathcal F}_2(z)$, autrement dit $\text{Im }W$ est tangent \`a $C_{\Delta(v)} \cap U^{\prime}$.  En particulier, en $\Delta(v)$, $\text{Im }W$ devra \^etre tangent aux orthogonaux de toutes les g\'en\'eratrices du c\^one de lumi\`ere. Ceci force $W_{\Delta(v)}=0$.  En reprenant cette information en compte dans la relation d'invariance de $W$, on obtient~:
 $$D_z\phi_X^{s_k}(W(u_{i,k}(z),u_{j,k}(z),u_{l,k}(z)))=o \left( \frac{1}{s_k} \right)$$
  si deux des indices $i,j,k$ valent $n$, sans \^etre tous les trois \'egaux \`a $n$, et 
  $$D_z\phi_X^{s_k}(W(u_{i,k}(z),u_{j,k}(z),u_{l,k}(z)))=o \left( \frac{1}{s_k^2} \right)$$ 
sinon. Le point $1(b)$ du Th\'eor\`eme \ref{thm.integrabilite-stable} conduit alors \`a  
$$W_z(u_i(z),u_j(z),u_l(z))=0,$$
 sauf \'eventuellement si deux des indices valent $n$, auquel cas $W_z(u_i(z),u_j(z),u_l(z))$ est colin\'eaire \`a $u_1(z)$.  Mais dans ce cas, les sym\'etries du tenseur de Weyl donnent :
$$g_z(W(u_n(z),u_j(z),u_n(z)),u_n(z))=g_z(W(u_n(z),u_n(z),u_n(z)),u_j(z))=0.$$ On obtient donc bien $W_z=0$.  \end{proof}

\subsection{Preuve du Th\'eor\`eme \ref{thm.semi-simple-analytique}}
\label{sec.preuve.semisimple}
Sous les hypoth\`eses du th\'eor\`eme,  le flot d'holono\-mie $\{ h^t \}$ de $X$ en $x_0$ poss\`ede une  partie lin\'eaire qui est semi-simple sur ${\bf C}$.  En particulier, quitte \`a conjuguer $\{ h^t \}$ dans $P$ par une translation, on peut supposer que $\{ h^t \}$ est le produit commutatif $\{ l^t \cdot \tau^t \}$ d'un flot lin\'eaire $\{ l^t \}$ semi-simple sur ${\bf C}$, et d'un flot de translations $\{ \tau^t \}$.  
Commen\c{c}ons par rappeler  un lemme de lin\'earisation, dont on peut par exemple trouver une preuve dans  \cite[Prop 4.2]{frances-riemannien}~:

\begin{lemme}
\label{lem.linearisable}
Soit $(M,g)$ une vari\'et\'e pseudo-riemannienne analytique, et $X$ un champ de vecteurs conforme analytique admettant une singularit\'e en $x_0$.  Si le flot d'holonomie $\{ h^t\}$ de $X$ en $x_0$ fixe un point  de $\BR^{p,q}$, alors $X$ et son champ d'holonomie $X_h$ sont analytiquement conjugu\'es au voisinage de $x_0$ et $o$ respectivement.
\end{lemme}

Si le flot de translations $\{ \tau^t \}$ est trivial, l'holonomie de $X$ fixe un point de ${\bf R}^{p,q}$.  Alors par ce lemme le champ $X$ est analytiquement conjugu\'e au champ de vecteurs lin\'eaire d\'efini par $\{ l^t \}$. Cela conduit au premier cas du Th\'eor\`eme \ref{thm.semi-simple-analytique}.  On suppose par la suite que $\{ \tau^t \}$ n'est pas trivial.

Soit $b_0$  un point de $B$ o\`u le flot d'holonomie de $X$ est de la forme $\{ l^t \cdot \tau^t \}$. L'hypoth\`ese d'analyticit\'e sur la vari\'et\'e $M$ implique que l'alg\`ebre de Lie des holonomies, en $b_0$, des champs conformes qui s'annulent en $x_0$ est alg\'ebrique. C'est une cons\'equence du Th\'eor\`eme de Frobenius pour les structures conformes \cite[Sec 3.4]{gromov.rgs} (voir aussi \cite[Thm 3.11]{me.frobenius}).  L'alg\`ebre de Lie d'un groupe alg\'ebrique est stable par la d\'ecomposition de Jordan, c'est-\`a-dire qu'elle contient les composantes semi-simples et nilpotentes de ses \'el\'ements (voir par exemple \cite{morris}, 4.4.2).  On en d\'eduit qu'il existe sur un voisinage $U$ de $x_0$, un champ conforme $Y$ dont le flot d'holonomie en $x_0$ est pr\'ecis\'ement $\{ \tau^t \}$.  On est alors dans le second cas du Th\'eor\`eme \ref{thm.translations}. Par cons\'equent, un ouvert non vide de $(M,g)$ va \^etre conform\'ement plat, et par analyticit\'e, $(M,g)$ est conform\'ement plate.  Le th\'eor\`eme d\'ecoule alors de la Remarque \ref{rem.cas-plat}.

%

\section{Le cadre lorentzien analytique : preuve du Th\'eor\`eme \ref{thm.lorentz-analytique}}
\label{sec.lorentz-analytique}

On consid\`ere un champ de vecteurs conforme analytique, sur une vari\'et\'e lorentzienne analytique $(M,g)$, et l'on suppose que le champ $X$ admet  une singularit\'e $x_0$.  La preuve du th\'eor\`eme va consister \`a analyser les holonomies possibles du champ $X$ en $x_0$.
Par le Lemme \ref{lem.linearisable}, si l'holonomie de $X$ en $x_0$ admet une d\'ecomposition affine dont la partie translation est triviale, alors $X$ est lin\'earisable au voisinage de $x_0$.  Si la partie lin\'eaire est compacte et la partie translation non triviale, alors la courbure conforme s'annule sur un ouvert non vide par le Th\'eor\`eme \ref{thm.translations}.  Il reste beaucoup  d'holonomies possibles entre ces deux situations.  Toutefois, pour les m\'etriques lorentziennes analytiques, on sera en mesure de  tra\^iter ces cas interm\'ediaires pour arriver au Th\'eor\`eme \ref{thm.lorentz-analytique}.  
 L'hypoth\`ese d'analyticit\'e va \^etre exploit\'ee, comme en  Section \ref{sec.preuve.semisimple}, pour dire que les parties semi-simples et unipotentes d'un flot d'holonomie, sont elles-m\^emes les flots d'holonomies de champs conformes au voisinage de $x_0$.  Cela va permettre de r\'eduire fortement le nombre de cas \`a \'etudier, soit par l'application du Th\'eor\`eme \ref{thm.translations}, qui assurera que $(M,g)$ est conform\'ement plate, soit en utilisant le Lemme \ref{lem.linearisable} pour conclure que $X$ est lin\'earisable.  Apr\`es ce travail pr\'eliminaire, il ne restera essentiellement que deux types d'holonomies \`a \'etudier.  Une analyse fine de la dynamique de $X$ au voisinage de $x_0$ dans ces deux cas permettra de prouver l'anulation de la courbure  conforme.  Ce sera l'objet des Sections \ref{sec.b.non.nul}  et \ref{sec.b.nul}.  

\subsection{R\'eduction de l'holonomie \`a deux cas}

Soit $\{ h^t \}$ le flot d'holonomie de $X$ en $x_0$, relativement \`a $b_0 \in B$ dans la fibre de $x_0$.  On consid\`ere la d\'ecomposition de Jordan de $\{h^t\}$ dans le groupe alg\'ebrique $P$~: $\{ h^t \}$ s'\'ecrit comme produit commutatif d'un flot $\{ h_s^t \}$ semi-simple sur $\BC$, c'est-\`a-dire que l'action de $\{ \Ad h_s^t \}$ sur $\liep$ est diagonalisable sur $\BC$, et d'un flot unipotent $\{ h_u^t \}$.  

Si le flot $\{ h_u^t \}$ est trivial, cela veut dire que $\{ h^t \}= \{ h_s^t \}$ est conjugu\'e dans $P$ \`a un flot de $\BR_+^* \times \OO(1,n-1)$.  L'holonomie fixe un point de $\BR^{1,n-1}$~: le Lemme \ref{lem.linearisable}  assure que $X$ est analytiquement lin\'earisable au voisinage de $x_0$.  

Supposons  maintenant  que $\{ h_u^t \}$ n'est pas trivial, mais qu'il fixe un point de $\BR^{1,n-1}$.  On appelle $E$ le sous-espace affine constitu\'e des points fixes de $\{ h_u^t \}$.  Le flot $\{ h_s^t \}$ agit sur $E$ comme un flot semi-simple dans le groupe $\Aff (E)$; il va donc avoir un point fixe, et finalement $\{ h^t \}$ va fixer un point de $\BR^{1,n-1}$~: on conclut \`a nouveau que $X$ est analytiquement lin\'earisable au voisinage de $x_0$.

Il reste \`a \'etudier le cas o\`u $\{ h_u^t \}$ ne fixe aucun point de $\BR^{1,n-1}$.  Comme nous l'avons expliqu\'e en \ref{sec.preuve.semisimple}, l'analycit\'e de $M$ permet de supposer que $\{ h^t \} = \{h_u^t \}$, ce que nous ferons par la suite.  On va raisonner au niveau de l'alg\`ebre de Lie $\liep$ et \'ecrire le champ d'holonomie $X_h$ comme une somme $X_h=U + T$, o\`u $U$ est un \'el\'ement nilpotent  de $\oo(1,n-1)$, et $T \in \lien^+$ est non trivial (sans quoi $h^t$ aurait un point fixe dans $\BR^{1,n-1}$).  Si $U$ est trivial, on conclut directement par le Th\'eor\`eme \ref{thm.translations} que $(M,g)$ est conform\'ement plate.




Si $U$ n'est pas trivial, alors,  \`a conjugaison dans $P$ pr\`es, on peut supposer que $U$ est la transformation suivante de $\BR^{1,n-1}$~:
$$ U : x \mapsto b(x,e_2) e_1 - b(x,e_n) e_2.$$ 
L'\'el\'ement  $T$ correspond \`a une  translation de vecteur $v$ sur $\BR^{1,n-1}$.

Soient $T_2$ et $T_n$ les \'el\'ements de $\lien^+$ correspondant aux  translations de $\BR^{1,n-1}$ de vecteurs  $e_2$ et $e_n$, respectivement.  Elles commutent avec  $T$.  Leurs commutateurs avec $U$ sont
$$ [T_2,U] = - T_1 \qquad \mbox{et} \qquad [T_n,U] = T_2$$
Quitte \`a conjuguer $X_h$ par une translation dans $P$ de direction dans $\mbox{Vect}(e_2,e_n)$, on peut alors supposer que $v$ n'a pas de composante dans $\mbox{Vect}(e_1,e_2)$.  \'Ecrivons $v = a \xi + be_n$ avec $\xi \in \mbox{Vect}(e_3, \ldots, e_{n-1})$ de norme $1$.  Quitte \`a conjuguer par une rotation de $P$ commutant avec l'exponentielle de $U$, on peut aussi supposer que $\xi = e_3$.

  Apr\`es ces diverses conjugaisons, l'\'el\'ement $X_h$ s'\'ecrit, dans $\oo(2,n)$, comme~: 
$$ 
X_h = 
\left(
\begin{array}{cccccc}
 0 & b & 0 & a \xi^*  &  &  \\
  &  0 & 1 &     &  &  \\
  &   & 0  &     & -1 & \\
  &   &   &  \ddots   &    & -a  \xi \\
  &   &   &     &   0 & -b  \\
  &   &   &     &    & 0
\end{array}
\right)
$$

En prenant l'exponentielle, on obtient~:
\begin{eqnarray} 
\label{eqn.ab-matrice}
h^t & = & 
\left(
\begin{array}{cccccc}
1 & tb & \frac{t^2b}{2} & ta\xi^*  & -\frac{t^3b}{6} & -\frac{t^2a^2}{2}  + \frac{b^2t^4}{24} \\
  & 1  & t &       & -\frac{t^2}{2} & \frac{t^3 b}{6} \\
  &    & 1 &       &      -t        & \frac{t^2 b}{2} \\
  &    &   &   \ddots    &                &  - ta\xi   \\
  &    &   &       &      1         & -tb   \\
  &    &   &       &                & 1
\end{array}
\right)
\end{eqnarray}          

Comme $T \neq 0$, alors $a$ et $b$ ne sont pas tous les deux nuls.

Nous allons \'etudier la dynamique de $\{ h^t \}$ pr\`es de $o$ sur $\Ein^{1,n-1}$, et utiliser les r\'esultats des sections pr\'ec\'edentes pour montrer que la courbure de Weyl s'annule sur un ouvert non vide de $(M,g)$.  L'analycit\'e donnera la platitude conforme globale de $M$.  Deux comportements dynamiques qualitativement diff\'erents apparaissent, suivant que $b$ est nul ou non.  Nous d\'etaillons \`a pr\'esent ces deux cas.

\subsection{Annulation de la courbure conforme dans le cas $b \neq 0$}
\label{sec.b.non.nul}
Quitte \`a conjuguer $\{ h^t \}$ dans $P$, nous supposerons ici que $b =1$ et $\xi=e_3$. 
Nous appelons  $\Omega:=\{ z \in  \Ein^{1,n-1} \ | \ z=j^o(x), \ q(x) \not = 0 \}$, et $U = B(o,R_0)$  le voisinage de $o$ donn\'e par la Proposition \ref{prop.stabilite}.
Le point cl\'e de cette section va \^etre de montrer la~:

\begin{proposition}
\label{prop.cas1}
Le flot $\{ h^t \}$ a les propri\'et\'es dynamiques suivantes~:
\begin{enumerate}
\item{Pour tout point $z \in \Omega$, on a $h^t.z \to o$ lorsque $t \to \pm \infty$, la convergence \'etant de plus uniforme sur les compacts de $\Omega$.  En particulier, pour toute suite de r\'eels $ (t_k)$ telle que $|t_k| \to \infty$, et pour tout $z \in \Omega$, la suite $(h^{t_k})$ est fortement stable en $z$.}
\item{Il existe un segment  g\'eod\'esique conforme $[\alpha]$ issu de $x_0$, se d\'eveloppant sur un segment g\'eod\'esique  $[\beta]$, avec $[\beta] \setminus \{ o \} \subset \Omega$, et tel que d'une part $h^t.[\beta] \subset U$ pour tout $t \geq 0$, et d'autre part $h^t.[\beta] \to o$ lorsque $t \to \infty$.}
\end{enumerate}

\end{proposition}
Cette proposition impliquera, par la Proposition \ref{prop.stabilite}, qu'il existe un ouvert $V$ de $(M,g)$ sur lequel le flot $\phi_X^t$ est d\'efini pour tout $t \geq 0$, ainsi qu'une suite $s_k \to \infty$ telle que $(\phi_X^{s_k})$ soit fortement stable sur $V$.  La Proposition \ref{prop.annulation} entra\^{\i}nera  la platitude conforme de $(M,g)$, prouvant ainsi le Th\'eor\`eme \ref{thm.lorentz-analytique}  dans ce cas.

\begin{proof} 
 Pour $x = (x_1, \ldots, x_n) \in \BR^{1,n-1}$, soit $\alpha_x$ le polyn\^ome de degr\'e $4$
$$ \alpha_x(t) = 1 + (x_1 + ax_3) \cdot t + \left( \frac{x_2}{2} + \frac{a^2 q(x)}{4} \right) \cdot t^2 - \frac{x_n}{6} \cdot t^3 - \frac{q(x)}{48} \cdot t^4$$

L'action du flot $\{ h^t \}$ sur $\Ein^{p,q}$ s'\'ecrit~:

\begin{eqnarray*}
\label{expression}
h^t.j^o(x) & = & \left[ 1 : \frac{1}{\alpha_x(t)} \cdot ( x_1 + x_2 t - \frac{x_n}{2}\cdot t^2 - \frac{q(x)}{12} \cdot t^3) : \right. \\
& &  \frac{1}{\alpha_x(t)} \cdot ( x_2 -  x_n t - \frac{q(x)}{4} \cdot t^2) : \frac{1}{\alpha_x(t)} \cdot (x_3 + \frac{a q(x)}{2} \cdot t) : \\
& &  \left. \frac{x_4}{\alpha_x(t)} : \cdots : \frac{x_{n-1}}{\alpha_x(t)} : \frac{1}{\alpha_x(t)} \cdot (x_n + \frac{q(x)}{2} \cdot t)  : - \frac{q(x)}{2 \alpha_x(t)} \right]
\end{eqnarray*}

Le premier point de la proposition d\'ecoule ais\'ement de cette expression.

Montrons \`a pr\'esent le second point.  Soit $z \in \Omega$.  D'apr\`es ce qui pr\'ec\`ede, il existe $T>0$ tel que si $t \geq T$, alors $h^t.z \in U$.  On d\'efinit $z_0=h^T.z$, $y_0 \in \BR^{1,n-1}$ par $j^o(y_0)=z_0$, et $y_t \in \BR^{1,n-1}$ par $h^t.j^o(y_0) = j^o(y_t)$ pour tout $t \geq 0$.
  On a bien entendu $h^t.z_0 \in U$ pour tout $t \geq 0$. Nous commen\c{c}ons par \'etablir une formule pour $h^t.(j^o(uy_0))$, $u \in \BR$.

\begin{lemme}
\label{lemme.formule-htux}
Pour $u \in \BR$,
$$ h^t.(j^o(uy_0)) = [ e_0 + u y_t - u \frac{q(y_t)}{2} e_{n+1} + u (1 - u) \frac{q(y_0)}{2} h^t.e_{n+1} ]$$
\end{lemme}

\begin{proof}
Rappelons la d\'efinition de $j^o$:
$$ j^o(y_t) = [ e_0 + y_t - \frac{q(y_t)}{2} e_{n+1} ]$$

Comme $h^t$ agit lin\'eairement  sur $\BR^{2,n}$, on a aussi
$$ h^t.j^o(y_0) = [ e_0 + h^t.y_0 - \frac{q(y_0)}{2} h^t.e_{n+1} ]$$

Donc
$$ h^t.y_0  = y_t - \frac{q(y_t)}{2} e_{n+1} + \frac{q(y_0)}{2} h^t.e_{n+1}$$

Alors pour $u \in \BR$,
\begin{eqnarray*}
h^t.(j^o(uy_0)) & = & [h^t.(e_0 + u y_0 - u^2\frac{q(y_0)}{2} e_{n+1} )]\\
 & = & [e_0 + u h^t.y_0 - u^2\frac{q(y_0)}{2} h^t.e_{n+1}] \\
 & = & [e_0 + u y_t - u \frac{q(y_t)}{2} e_{n+1} + u (1 - u) \frac{q(y_0)}{2} h^t.e_{n+1}]
 \end{eqnarray*} \end{proof}

   Nous d\'efinissons la g\'eod\'esique conforme $\gamma : [0,1] \to \Ein^{1,n-1}$ par $\gamma(s):=j^o(-sy_0)$.  Nous allons montrer que $h^t.[\gamma] \to o$ lorsque $t \to \infty$.
Si tel n'\'etait pas le cas,  il existerait une  suite $\{u_k\}$  de  $[-1,0)$, un voisinage $W$ de $o$ et $t_k \to \infty$ tels que pour tout $k \in \BN$, $h^{t_k}.(j^o(u_k y_0)) \notin W$.     Notons que la suite $(u_k)$ tend vers $0$, car sinon, on aurait une contradiction avec le premier point de la proposition.     Quitte \`a consid\'erer une suite extraite, nous pouvons supposer que $(u_kt_k^4)$ admet une limite dans $[-\infty,0]$.  Si cette limite est dans $(- \infty,0]$, alors on tire de (\ref{eqn.ab-matrice}) que $u_k h^{t_k}.e_{n+1}$ tends vers $c e_0$ pour $c \leq 0$ quand $k \rightarrow \infty$.   Par ailleurs, $h^t.j^o(y_0) \rightarrow o$, et donc
 $$ \lim_{t \rightarrow \infty} y_t = \lim_{t \rightarrow \infty} q(y_t) = 0.$$
 
 On d\'eduit du Lemme \ref{lemme.formule-htux} que $h^{t_k}.j^o(u_ky_0)$ tend vers $[(1+c\frac{q(y_0)}{2})e_0]=o$ quand $k \rightarrow \infty$, en contradiction avec l'hypoth\`ese $h^{t_k}.(j^o(u_k x)) \notin W$.   Une remarque cl\'e est que $q(y_0)<0$, et donc,  $(1+c\frac{q(y_0)}{2})e_0 \not = 0$.  En effet, $q(y_0) > 0$ impliquerait  $\alpha_{y_0}(t) < 0$ pour $t$ suffisamment grand car le coefficient de $t^4$ dans $\alpha_{y_0}(t)$ est $- q(y_0)/48$.  Mais $\alpha_{y_0}(0) = 1$.  Donc $\alpha_{y_0}$ s'annulerait pour un $t_0 \in [0,\infty)$, ce qui contredirait $h^{t_0}.z \in U$.   
 
 Il ne reste plus qu'\`a consid\'erer le cas o\`u $u_kt_k^4 \rightarrow - \infty$.   Comme les composantes de $h^{t_k}.(j^o(u_k y_0))$, sauf celle sur $e_0$,  croissent au plus  comme $u_k t_k^3$, elles sont n\'egligeables devant la composante selon $e_0$ qui, elle, est de l'ordre de $u_kt_k^4$~: la limite est encore $o$, d'o\`u une nouvelle contradiction.  
 
Du fait que $h^t.[\gamma] \to o$, on d\'eduit l'existence de $T^{\prime}>0$ tel que $h^t.[\gamma] \subset U$ pour tout $t \geq T^{\prime}$.  D\'efinissons $\beta(s):=h^{T^{\prime}}.\gamma(s)$, $s \in [0,1]$.  On peut  \'ecrire $\beta(s)=\pi_G(e^{-s  \xi})$, o\`u $\xi \in \lieg$. Quitte \`a remplacer $s \mapsto \beta(s)$ par $s \mapsto \beta(\delta s)$, avec $\delta >0$, on a que $\beta$ est le d\'eveloppement de la g\'eod\'esique conforme  $\alpha(s):=\pi(\exp(b_0,-s \xi))$, d\'efinie sur $[0,1]$.  Le second point de la Proposition \ref{prop.cas1} est prouv\'e.  \end{proof}

\subsection{Annulation de la courbure de Weyl dans le cas $b=0$}
\label{sec.b.nul}

Il s'agit du cas le plus difficile.  On va montrer gr\^ace \`a la Proposition \ref{prop.stabilite} que le flot $\{ \phi_X^t \}$ est stable, mais pas fortement stable, sur un ouvert $V$, et que cet ouvert est ``\'ecras\'e" sur un segment g\'eod\'esique de lumi\`ere, fix\'e par le flot (voir le Lemme \ref{lem.dynamique.phi}).  Ceci n'est pas suffisant a priori pour montrer que $V$ est conform\'ement plat.  Toutefois, une analyse plus fine de la dynamique,  et l'emploi du Th\'eor\`eme \ref{thm.integrabilite-stable} vont nous permettre  de montrer l'annulation du tenseur de Weyl sur un ouvert.  La preuve va requ\'erir quatre \'etapes que nous d\'etaillons ci-dessous.

Puisque $b = 0$, alors $a$ est forc\'ement non nul, car nous avons exclu le cas o\`u $a$ et $b$ valent simultan\'ement $0$.  On peut encore supposer, quitte \`a conjuguer $\{h^t\}$ par un \'el\'ement du centralisateur de $\{ e^{tU} \}$ dans $P$, que $a = 1$ et $\xi = e_3$.  Le flot $\{ h^t \}$ admet alors l'expression matricielle~:

\begin{equation}
\label{eq.expression.matricielle}
 h^t   =  
\left(
\begin{array}{cccccc}
1 & 0 & 0 & t \xi^* & 0 & -\frac{t^2}{2}   \\
  & 1  & t &       & -\frac{t^2}{2} & 0 \\
  &    & 1 &       &      -t        & 0  \\
  &    &   &   \ddots    &                &  - t \xi  \\
  &    &   &       &      1         & 0  \\
  &    &   &       &                & 1
\end{array}
\right)
\end{equation}

Nous voyons que le flot $\{ h^t \}$  fixe point par point la g\'eod\'esique de lumi\`ere $\Lambda$, projectivis\'ee dans $\Ein^{1,n-1}$ du plan ${\text Vect} (e_0,e_1) \subset \BR^{2,n}$.  Cette g\'eod\'esique de lumi\`ere peut \^etre param\'etr\'ee au voisinage de $o$ par $\Lambda(s):= \pi_G(e^{s \xi_1})$, o\`u $\xi_1 \in \lien^{-}$.


\subsubsection{Premi\`ere \'etape : dynamique de l'holonomie dans le mod\`ele}
Commen\c{c}ons tout d'abord par quelques consid\'erations g\'eom\'etriques.  
 Nous  d\'esi\-gnons par $\Omega_{\Lambda}$ l'ouvert $\Ein^{1,n-1} \setminus \Lambda$.  Il est facile de v\'erifier que si $z \in \Omega_{\Lambda}$, alors $C(z)$, le c\^one de lumi\`ere issu de $z$,  coupe $\Lambda$ en un unique point, que l'on note $\pi_{\Lambda}(z)$.  On h\'erite ainsi d'une submersion $\pi_{\Lambda} : \Omega_{\Lambda} \to \Lambda$, dont les fibres sont des hypersurfaces d\'eg\'en\'er\'ees (ces fibres sont les intersections des c\^ones de lumi\`ere de la forme $C(x)$, $x \in \Lambda$, avec $\Omega_{\Lambda}$).  On appelera $F_{\Lambda}$ le feuilletage de $\Omega_{\Lambda}$ par les fibres de l'application $\pi_{\Lambda}$.

Soit $\tau$ la transformation de $\Ein^{1,n-1}$ d\'efinie par $[x_0,x_1,x_2, \ldots,x_{n+1}] \mapsto [x_1,-x_0,x_2, \ldots ,x_{n+1}].$  C'est une application conforme qui laisse $\Lambda$ globalement invariante, et qui agit sans point fixe sur $\Lambda$.
Dans la suite, on appelle $U := B(o,R_0)$ le voisinage de $o$ donn\'e par la Proposition \ref{prop.stabilite}.

\begin{proposition}
\label{prop.dynamique.2}
L'action de $\{ h^t \}$ sur $\Ein^{1,n-1}$ a les propri\'et\'es suivantes~:
\begin{enumerate}

\item{Pour tout $z \in \Omega_{\Lambda}$, $\lim_{t \to \pm \infty} h^t.z=\tau(\pi_{\Lambda}(z))$, la convergence \'etant uniforme sur les compacts de $\Omega_{\Lambda}$.}
\item{Pour tout $z \in \Omega_{\Lambda}$, et toute suite $(t_k )$ telle que $|t_k| \to \infty$,  $(h^{t_k})$ est stable en $z$, mais pas fortement stable.}
\item{Il existe un segment  g\'eod\'esique conforme $[\alpha]$ issu de $x_0$, se d\'eveloppant sur un segment g\'eod\'esique  $[\beta]$, avec $[\beta] \setminus \{ o \} \subset \Omega_{\Lambda}$, et tel que d'une part $h^t.[\beta] \subset U$ pour tout $t \geq 0$, et d'autre part $h^t.[\beta] \to o$ lorsque $t \to \infty$.}
\end{enumerate}
\end{proposition}

\begin{proof}
Consid\'erons $z \in \Omega_{\Lambda}$.  On \'ecrit $z:=[x_0: \cdots : x_{n+1}]$, avec $x_n$ et $x_{n+1}$ qui ne sont pas tous les deux nuls.  De l'expression matricielle (\ref{eq.expression.matricielle}), on tire ais\'ement que 
$$\lim_{t \to \pm \infty} h^t.z=[x_{n+1}:x_n:0:\cdots : 0].$$
La convergence est de plus uniforme sur les compacts de $\Omega_{\Lambda}$.  
Or $\pi_{\Lambda}(z)=[x_n:- x_{n+1}:0: \cdots :0].$  On a donc bien $\lim_{t \to \pm \infty}=\tau(\pi_{\Lambda}(z)).$ Ceci montre le premier point de la proposition.  La convergence uniforme de $h^t$ vers $\tau \circ \pi_{\Lambda}$ sur les compacts de $\Omega_{\Lambda}$ montre que pour toute suite $(t_k)$ telle que $|t_k| \to \infty$, et tout $z \in \Omega_{\Lambda}$, la suite $( h^{t_k} )$ est stable en $z$.  Elle n'est pas fortement stable car les fibres de l'application $\tau \circ \pi_{\Lambda}$ sont les feuilles du feuilletage ${F}_{\Lambda}$.  En particulier, il existe des points $z^{\prime}$ arbitrairement proches de $z$ tels que $\tau(\pi_{\Lambda}(z^{\prime})) \not =\tau(\pi_{\Lambda}(z))$.  Il n'existe donc pas de voisinage $V$ de $z$ tel que  $h^{t_k}(V) \to \tau(\pi_{\Lambda}(z)).$

Il nous reste \`a montrer le dernier point de la proposition.
Pour $x \in \BR^{1,n-1}$ fix\'e, on introduit le polyn\^ome $\alpha_x$ suivant
$$ \alpha_x(t) = 1 + x_3 t+ \frac{q(x)}{4} \cdot t^2.$$

Alors le flot de $\{ h^t \}$ sur $\Ein^{1,n-1}$ s'\'ecrit
\begin{eqnarray*}
h^t.j^o(x) & = & \left[ 1 : \frac{1}{\alpha_x(t)} \cdot ( x_1 + x_2 t - \frac{x_n}{2}\cdot t^2 ) :   \frac{1}{\alpha_x(t)} \cdot ( x_2 -  x_n t  )  \right. \\
&  & \left.  : \frac{1}{\alpha_x(t)} \cdot (x_3 + \frac{q(x)}{2} \cdot t) :  \frac{x_4}{\alpha_x(t)} : \cdots : \frac{x_n}{\alpha_x(t)}  : - \frac{q(x)}{2 \alpha_x(t)} \right]
\end{eqnarray*}

On choisit $x \in \BR^{1,n-1}$ tel que   $q(x) \neq 0$ et $x_n=0$. Alors $h^t.j^o(x) \rightarrow o$, donc il existe $T>0$ tel que pour $t \geq T$, on a $h^t.j^o(x) \in U$. On d\'efinit par la suite $y_0$ par $j^o(y_0)=h^T.j^o(x)$.  On a $h^t.j^o(y_0) \in U$ pour tout $t \geq 0$, et donc $\alpha_{y_0}(t) \not =0$ pour $t \geq 0$.  Il s'ensuit que $q(y_0)>0$.   

Par le Lemme \ref{lemme.formule-htux}, pour $u \in \BR$,
$$ h^t.(j^o(uy_0)) = [e_0 + u y_t - u \frac{q(y_t)}{2} e_{n+1} + u (1 - u) \frac{q(y_0)}{2} h^t.e_{n+1}]$$
o\`u l'on a \`a nouveau pos\'e $h^t.j^o(y_0) = j^o(y_t)$.  On a encore  $y_t \rightarrow 0$, et donc $q(y_t) \to 0$.  

Posons $\gamma(s):=j^o(-sy_0)$ pour $s \in [0,1]$. On remarque que $[\gamma] \setminus \{ o \} \subset \Omega_{\Lambda}$. Nous allons montrer que $h^t.[\gamma] \to o$.  On concluera alors exactement de la m\^eme mani\`ere que pour la fin de la Proposition \ref{prop.cas1}.

Pour montrer que  $h^t.[\gamma] \to o$, il suffit de montrer que pour toute suite $( u_k )$ de $[-1,0)$, et toute suite $t_k \to \infty$, on a, quitte \`a consid\'erer une sous-suite,  $\lim_{k \to \infty}h^{t_k}.\gamma(u_ky_0)=o$.  

Si la suite $(u_kt_k^2)$ est born\'ee, on prend une sous-suite de sorte qu'elle converge vers  $l \in (- \infty,0]$. Dans ce cas $1-\frac{lq(y_0)}{4} > 0$ et on obtient, en utilisant l'expression (\ref{eq.expression.matricielle}), que $h^{t_k}.j^o(u_k y_0) \rightarrow o$.  

Si $( u_k t_k^2 )$ n'est pas born\'ee, on  peut supposer qu'elle tend vers $- \infty$.   Toutes les composantes de $h^{t_k}.j^o(u_ky_0)$ croissent en $O(u_kt_k)$, sauf celle selon $e_0$ qui est de l'ordre de $u_k t_k^2$~: la limite  de $(h^{t_k}.j^o(u_ky_0))$ est encore $o$.   \end{proof}


\begin{corollaire}
\label{coro.holonomie.h}
Pour toute suite $( t_k )$ telle que $|t_k| \to \infty$, et tout $z \in \Omega_{\Lambda}$, la suite $( h^{t_k} )$ admet, en $z$, une holonomie de la forme\\
 $$(h_k)=(\text{diag} (\lambda_2(k)^2, \lambda_2(k), \ldots, \lambda_2(k), 1)), \ \text{o\`u} \ 1/\lambda_2(k) \to 0.$$
\end{corollaire}

\begin{proof}
On sait par le second point de la Proposition \ref{prop.dynamique.2} que $(h^{t_k} )$ est stable en $z$.  Comme nous l'avons d\'ej\`a mentionn\'e, cela signifie qu'elle admet en $z$ une suite d'holonomie de la forme $(\diag(\lambda_1(k), \ldots , \lambda_n(k))) \in (A^+)^{\BN}$ (voir \cite[Lemme 4.3]{frances-degenere}).  Mais nous sommes ici en signature lorentzienne, ce qui veut dire qu'il existe $\sigma_k$ et $\mu_k$ positifs, $\sigma_k \geq \mu_k \geq 1$ tels que $\lambda_1(k)=\sigma_k\mu_k$, $\lambda_i(k)=\mu_k$ si $i =2, \ldots ,n-1$ et $\lambda_n(k)=\frac{\mu_k}{\sigma_k}$. Par ailleurs, toujours par le second point de  la Proposition \ref{prop.dynamique.2}, on sait que $( h^{t_k} )$ n'est pas fortement stable en $z$, et il en va de m\^eme pour toutes ses sous-suites.  Aussi, $1/\lambda_n(k) \geq \delta>0$ pour tout $k$.  Quitte \`a multiplier l'holonomie par une suite born\'ee de $P$ et prendre une sous-suite, on peut donc supposer que $\lambda_n(k)=1$.   On obtient alors une holonomie $(\diag (\lambda_2(k)^2, \lambda_2(k), \ldots, \lambda_2(k), 1))$ de la forme annonc\'ee. \end{proof}

\subsubsection{Deuxi\`eme \'etape: propri\'et\'es dynamiques de $\{\phi_X^t\}$ au voisinage de $\alpha(1)$ }

Nous reprenons les conclusions et les notations de la Proposition \ref{prop.dynamique.2}.  En particulier, dans ce qui suit, le segment g\'eod\'esique $[\alpha]$ est celui donn\'e par le troisi\`eme point de  la Proposition \ref{prop.dynamique.2}. 
\begin{lemme}
\label{lem.dynamique.phi}
Soit $x_1=\alpha(1)$.  Il existe un voisinage $V$ de $x_1$, relativement compact dans $M$,  sur lequel $\phi_X^t$ est d\'efini pour tout $t \geq 0$.  De plus, il existe une suite $( s_k )$ qui tend vers l'infini, et une submersion  lisse $\rho : V \rightarrow \Delta$ telles que~:
\begin{enumerate}
\item{Pour tout $x \in V$, $\phi^{s_k}_X.x \rightarrow \rho(x)$.}
\item{La suite $( \phi_X^{s_k} )$ est stable en chaque point de $V$, et admet une holonomie de la forme $(h_k) =(\diag (\lambda_2(k)^2, \lambda_2(k), \ldots, \lambda_2(k),1))$ o\`u $\frac{1}{\lambda_2(k)} \to 0$.}
\end{enumerate}

\end{lemme}

\begin{proof}
Par la Proposition \ref{prop.dynamique.2}, $[\alpha]$ satisfait aux hypoth\`eses de la Proposition \ref{prop.stabilite}.  On peut donc affirmer qu'il existe un voisinage $V$ de $x_1:=\alpha(1)$ tel que $\phi_X^t$ soit d\'efini sur $V$ pour tout $t \geq 0$.   Par ailleurs $\lim_{t \to \infty} \phi_X^t.[\alpha]=x_0$, et il existe une suite $s_k \to \infty$  telle que la suite d'holonomie de $( \phi_X^{s_k}) $ en chaque point de $V$  soit l'holonomie de $(h^{s_k})$ en $\beta(1)$.  Par le Corollaire \ref{coro.holonomie.h}, cette holonomie est de la forme $(h_k)=(\diag(\lambda_2(k)^2, \lambda_2(k), \ldots, \lambda_2(k), 1))$ avec $\frac{1}{\lambda_2(k)} \to 0$, ce qui prouve le point (2).  

Si l'on consid\`ere $(h_k)$ comme une suite de $\OO(p+1,q+1)$, alors $(\Ad h_k)$ restreinte \`a $\mathfrak{n}^-$ est equivalente \`a l'action de $(h_k)$ sur $\BR^{1,n-1}$ par la representation standard transpos\'ee.  Cette suite adjointe tend vers une application, not\'ee $L_{\infty} \in \End (\lien^-)$, avec $\text{Im } L_{\infty}=\BR. \xi_1$  (o\`u $\xi_1$ est le vecteur de $\lien^-$ tel que $\Lambda(s)=\pi_G(e^{s \xi_1})$).  On peut restreindre $V$ de sorte qu'il soit relativement compact.

Choisissons $b_1 \in B$ au-dessus de $x_1$.  Quitte \`a restreindre encore $V$, on peut supposer que $V=\pi(\exp(b_1,{\mathcal V}))$, o\`u ${\mathcal V}$ est un voisinage relativement compact de $0$ dans $\lien^-$  et $\zeta \mapsto \pi (\exp(b_1, \zeta))$ r\'ealise un diff\'eomorphisme de $\mathcal V$ sur $V$.  Comme $( h_k )$ est une suite d'holonomie en $x_1$, et comme $\phi_X^{s_k}. [\alpha] \to x_0$, il existe une suite $(b_k )$ qui converge vers $b_1$ telle que $b_k^{\prime}:=\phi_X^{s_k}.b_k.h_k^{-1}$ converge,  lorsque $k \to \infty$, vers $b_0^{\prime}$ dans la fibre de $x_0$.  Si $x \in V$, on \'ecrit $x=\pi(\exp(b_1,\zeta))=\pi(\exp(b_k,\zeta_k))$, o\`u $\zeta \in {\mathcal V}$, et $\zeta_k \rightarrow \zeta$.  On obtient alors pour tout $k$, en utilisant la relation (\ref{eq.f.exponentielle}) donn\'ee en la Section \ref{sec.transformations}~:
$$ \phi_X^{s_k}.\exp(b_k, \zeta_k).h_k^{-1}=\exp(b_k^{\prime}, (\Ad h_k)(\zeta_k)),$$ 
ou encore, en projetant sur $M$~:
$$ \lim_{k \to \infty} \phi_X^{s_k}(\pi(\exp(b_1,\zeta)))=\pi(\exp(b_0^{\prime}, L_{\infty}(\zeta))).$$
Appelons $\tilde{\Delta}(s):=\pi(\exp(b_0^{\prime}, s \xi_1))$, pour $s \in (-\delta,\delta)$, avec $\delta>0$ assez petit.  Quitte \`a restreindre encore $V$, l'application 
$$\rho: \pi(\exp(b_1,\zeta)) \mapsto \pi(\exp(b_0^{\prime}, L_{\infty}(\zeta)))$$ est une submersion de $V$ sur un intervalle de $\tilde \Delta$.  

Pour terminer la preuve du lemme, il ne nous reste plus qu'\`a montrer que $\tilde \Delta=\Delta$ (en tant que segments g\'eod\'esiques de lumi\`ere, abstraction faite du param\'etrage).  
Si  $\epsilon>0$ est suffisamment petit, alors il existe un petit voisinage $V^{\prime}$ de $x_1$, tel que $\phi^t_X.x \in V$ pour tout $x \in V^{\prime}$ et $t \in (-\epsilon, \epsilon)$.  Soit $I = \rho(V)$ et $I' = \rho(V')$; ce sont deux ouverts de $\tilde{\Delta}$ qui contiennent $x_0$.  De l'identit\'e $\lim_{k \to \infty} \phi_X^{s_k}.\phi_X^t.x=\phi_X^t.(\lim_{k \to \infty} \phi_X^{s_k}.x)$ pour tout $x \in V^{\prime}$, on d\'eduit que $\phi_X^t.I' \subseteq I$ pour tout $t \in (-\epsilon, \epsilon)$.  On conclut que $\phi_X^t.{\tilde \Delta} \subset {\tilde \Delta}$, et en particulier, $D_{x_0}\phi_X^t({\tilde \Delta}^{\prime}(0)) \in \BR {\tilde \Delta}^{\prime}(0)$.  Si l'on \'ecrit, en utilisant les relations (\ref{eq.iota}) et (\ref{eq.differentielle.f}) de \ref{sec.transformations}, que ${\tilde \Delta}^{\prime}(0)=\iota_{b_0}({\overline \zeta}_1)$, on obtient 
$\iota_{b_0}((\overline{\mbox{Ad}}\ h^t).(\overline{\zeta}_1)) \in \BR \iota_{b_0}(\overline{\zeta}_1)$, 
et ce  pour tout $t \in (-\epsilon, \epsilon)$.  Appelons $\overline{\lambda}^{1,n-1}$ la m\'etrique lorentzienne induite par $\lambda^{1,n-1}$ sur $\lieg/\liep$  (voir la Section \ref{sec.convexe} pour ces notations). Il est facile de v\'erifier que $\BR \overline{\xi}_1$ est la seule direction  de $\lieg/\liep$ qui soit isotrope relativement \`a  $\overline{\lambda}^{1,n-1}$, et  invariante par $\overline{\mbox{Ad}} \ h^t$.  On obtient donc $\overline{\zeta}_1=\overline{\xi}_1$, et  $\tilde{\Delta}^{\prime}(0)=\Delta^{\prime}(0)$.  
Remarquons pour conclure qu'en g\'eom\'etrie pseudo-riemannienne, deux segments g\'eod\'esiques de lumi\`ere qui ont une m\^eme tangente en un point sont identiques sur l'intersection de leurs domaines de d\'efinition.   \end{proof}

Notons que lorsque la dimension de $M$ est $3$, on d\'eduit directement du second point du lemme ci-dessus que $V$ est un ouvert conform\'ement plat, et par analyticit\'e, $(M,g)$ est aussi conform\'ement plate.  Cela r\'esulte de la Proposition 5 de \cite{frances-causal}.  Nous supposerons dor\'enavant que la dimension de $M$ est au moins quatre.

Dans ce qui suit, nous noterons $\PL(TM)$ le projectivis\'e du fibr\'e des vecteurs non nuls de type lumi\`ere de $TM$.  On d\'efinira $\PL(V)$ de mani\`ere pareille lorsque $V$ est un espace vectoriel muni d'un produit scalaire lorentzien. Dans la proposition qui suit, $\alpha$ et $\beta$ sont les  g\'eod\'esiques conformes donn\'ees par le troisi\`eme point de la Proposition \ref{prop.dynamique.2}.

\begin{proposition}
\label{prop.limite.stable}
Si $x_1:=\alpha(1)$, et $[u] \in \PL(T_{x_1}M)$, on a~:
$$\lim_{k \to \infty} D_{x_1} \phi_X^{s_k}([u])=[\Delta^{\prime}(0)],$$
la limite \'etant prise dans $\PL(TM)$.

\end{proposition}

\begin{proof}
\'Ecrivons $\alpha(v)=\pi(\exp(b_0,v\xi_{\alpha}))$, o\`u $\xi_{\alpha} \in (\mbox{Ad} \ P)(\lien^-)$.  On pose $\hat \alpha (v)= \exp(b_0,v\xi_{\alpha})$.  Le d\'eveloppement de $\hat \alpha$ en $1_G$ est $\hat \beta$ d\'efinie par $\hat \beta (v)=e^{v \xi_{\alpha}}$.  Enfin $\beta=\pi_G \circ \hat \beta$. Nous savons que $h^{s_k}.[\beta] \to o$.  On en d\'eduit l'existence d'une courbe $p_k : [0,1] \to P$ telle que $h^{s_k}{\hat{\beta}}(v)h^{-s_k}p_k(v)^{-1}$ soit une courbe du groupe $N^-:=e^{\lien^-}$.  Comme $\pi_G \circ e$ r\'ealise un diff\'eomorphisme de $\lien^-$ sur un voisinage de $o$ dans $\Ein^{1,n-1}$, on obtient~:
$$ h^{s_k}\hat{\beta}(1)(p_k h^{s_k})^{-1} \to 1_G,$$
o\`u l'on a pos\'e $p_k=p_k(1)$.
On en d\'eduit que $( \tilde{h}_k) :=(p_k h^{s_k})$ est une suite d'holonomie de $(h^{s_k})$ au point $\hat{\beta}(1)$.

Dans ce qui suit, l'\'el\'ement  de $\PL (\lieg/\liep)$ d\'etermin\'e par un $\overline{\zeta}$ non nul de $\lieg/\liep$  sera not\'e $[\zeta]$.   Nous allons montrer~:

\begin{lemme}
\label{lem.action.projectivise}
Pour tout $[\zeta] \in \PL (\lieg/\liep)$, on a $\lim_{k \to \infty} (\overline{\mbox{Ad}}\ \tilde{h}_k).[\zeta] \to [\xi_1]$.
\end{lemme}

\begin{proof}
Dans l'univers d'Einstein $\Ein^{1,n-1}$, le projectivis\'e du fibr\'e des vecteurs de type lumi\`ere, not\'e $\PL (T \Ein^{1,n-1})$, s'identifie \`a l'espace des g\'eod\'esi\-ques de lumi\`eres marqu\'ees de $\Ein^{1,n-1}$.  Maintenant, si $\Gamma$ est une g\'eod\'esique de lumi\`ere passant par $\beta(1)$, on aura que $h^{s_k}.\Gamma \to \Lambda$.  Pour le voir, commen\c{c}ons par supposer que  $\Gamma$ ne passe pas par $\pi_{\Lambda}(\beta(1))$.  Alors pour $z \in \Gamma$ diff\'erent de $\beta(1)$, $\pi_{\Lambda}(z) \not = \pi_{\Lambda}(\beta(1))$.  Il ressort de la Proposition \ref{prop.dynamique.2} que $\lim_{k \to \infty} h^{s_k}.\beta(1) \not = \lim_{k \to \infty} h^{s_k}.z$, o\`u ces deux limites sont  des points de $\Lambda$.  On conclut donc que $\lim_{{k \to \infty}}h^{s_k}.\Gamma=\Lambda$.
Maintenant, si $\Gamma$ passe par $\beta(1)$ et $\pi_{\Lambda}(\beta(1))$, on a d'une part $\lim_{k \to \infty} h^{s_k}.\beta(1)=\tau(\pi_{\Lambda}(\beta(1)))$, et d'autre part $\lim_{k \to \infty} h^{s_k}.\pi_{\Lambda}(\beta(1))=\pi_{\Lambda}(\beta(1))$, puisque $\pi_{\Lambda}(\beta(1)) \in \Lambda$.  Comme $\tau(\pi_{\Lambda}(\beta(1))) \not = \pi_{\Lambda}(\beta(1))$, on a l\`a encore $\lim_{{k \to \infty}}h^{s_k}.\Gamma=\Lambda$.  Si l'on exprime ceci dans $\PL (T \Ein^{1,n-1})$, cela signifie que pour tout $[u] \in \PL(T_{\beta(1)} \Ein^{1,n-1})$~:
\begin{equation}
\label{eq.limite.einstein}
 \lim_{k \to \infty} D_{\beta(1)} h^{s_k}([u])=[\Lambda^{\prime}(0)].
 \end{equation}
\'Ecrivons $[u]=\iota_{\hat{\beta}(1)}([\zeta])$.  Alors, utilisant les relations (\ref{eq.iota}) et (\ref{eq.differentielle.f})~:
$$D_{\beta(1)}h^{s_k}(\iota_{\hat{\beta(1)}}([\zeta]))=\iota_{h^{s_k}\hat{\beta}(1){\tilde h}_k^{-1}}((\overline{\mbox{Ad}}\ \tilde{h}_k).[\zeta]).$$
Comme nous avons vu que $(h^{s_k}\hat{\beta}(1){\tilde h}_k^{-1})$ tend vers $1_G$, la relation (\ref{eq.limite.einstein}) va impliquer $(\overline{\mbox{Ad}} \ \tilde{h}_k).[\zeta] \to [\xi_1]$.   \end{proof}

\`A $k \in \BN$ fix\'e, on consid\`ere la courbe~
$$v \mapsto \phi_X^{s_k}.\exp(b_0, v \xi_{\alpha}).h_k^{-1}=\exp(\phi_X^{s_k}.b_0.h_k^{-1}, (\mbox{Ad } h_k)( v \xi_{\alpha})).$$
  Le d\'eveloppement de cette courbe en $1_G$ est  $v \mapsto h^{s_k}e^{ v \xi_{\alpha}}h^{-s_k}$.  Par la formule (\ref{equ.formule-Cartan}) donn\'ee en Section \ref{sec.exponentielle},  la courbe $\hat \gamma_k: v \mapsto  \phi_X^{s_k}.\hat{\alpha}(v). h^{-s_k}p_k(v)^{-1}$ se d\'eveloppe en $1_G$ sur $\hat \sigma_k: v \mapsto h^{s_k}{\hat{\beta}}(v).h^{-s_k}p_k(v)^{-1}$.   Donnons nous sur le groupe $G$ une m\'etrique riemannienne invariante \`a gauche, d\'efinie \`a partir d'un produit scalaire $<\ , \ >_{\lieg}$ sur $\lieg$.  On appelle $L^G(\hat \sigma_k)$ la longueur de $\hat{\sigma}_k$ pour cette m\'etrique.  Sur $B$, on d\'efinit une m\'etrique riemannienne $< \ , \ >_B$ en posant $<u,w>_{B}=<\omega(u),\omega(w)>_{\lieg}$.  On appelle $L(\hat \gamma_k)$ la longueur de $\hat \gamma_k$ relativement \`a cette m\'etrique.  Il est clair que $L(\hat \gamma_k)=L^G(\hat \sigma_k)$. Comme $\hat \sigma_k$ est une courbe du groupe $N^-$, qui se projette dans $\Ein^{1,n-1}$ vers $h^{s_k}.[\beta]$, qui tend vers $o$, il d\'ecoule de \cite[Prop 3.2]{frances-riemannien} que $L^G(\hat \sigma_k)$ tend vers $0$, et donc $[\hat \gamma_k]$ tend vers $b_0$. 
 En particulier~:
$$ \lim_{k \to \infty} \phi_X^{s_k}.\hat \alpha(1). {\tilde h}_k^{-1}=b_0.$$

Pour tout $[\zeta] \in \PL(\lieg/\liep)$, on peut alors \'ecrire~:
$$D_{x_1}\phi_X^{s_k}(\iota_{\hat{\alpha}(1)}([\zeta]))=\iota_{\phi_X^{s_k}.\hat{\alpha}(1).{\tilde h}_k^{-1}}((\overline{\mbox{Ad}} \ \tilde{h}_k).[\zeta]),$$
ce qui, au vu du Lemme \ref{lem.action.projectivise}, conduit, pour tout $[u] \in \PL(T_{x_1}M)$, \`a~:
$$ \lim_{k \to \infty} D_{x_1}\phi_X^{s_k}([u])=\iota_{b_0}([\xi_1])=[\Delta^{\prime}(0)].$$
Ceci conclut la preuve de la proposition. \end{proof}

\subsubsection{Troisi\`eme \'etape : annulation du tenseur de  Weyl en  $x_0$}
\label{sec.annulation.x0}
Nous repre\-nons les conclusions du Lemme \ref{lem.dynamique.phi}~: il existe un ouvert $V$ contenant $x_1=\alpha(1)$ sur lequel $\phi_X^t$ est d\'efini pour $t \geq 0$, et une suite $s_k \to \infty$ telle que l'holonomie de $( \phi_X^{s_k} )$ en tout point de $V$ soit de la forme $$(h_k)= (\diag(\lambda_2(k)^2,\lambda_2(k), \ldots , \lambda_2(k), 1)), \text{ avec } 1/\lambda_2(k) \to 0.$$  On va appliquer le Th\'eor\`eme \ref{thm.integrabilite-stable} \`a $( \phi_X^{s_k} )$, quitte \`a restreindre l'ouvert $V$.  L'entier $s$ donn\'e par le Th\'eor\`eme \ref{thm.integrabilite-stable} vaut ici $3$.  Les suites $(\mu_i(k))$ fournies par le th\'eor\`eme sont $\mu_3(k)=1$, $\mu_2(k) = 1/\lambda_2(k)$, et $\mu_1(k)=\mu_2(k)^2$ pour tout $k$.  De la Remarque \ref{rem.orthogonalite}, on tire que le feuilletage $F_2$ de $V$ est un feuilletage par hypersurfaces d\'eg\'en\'er\'ees, et que $F_1$ est un feuilletage par g\'eod\'esiques de lumi\`eres. Pour tout $x$ dans $V$, on a de plus ${\mathcal F}_1(x)^{\perp}={\mathcal F}_2(x)$.

On consid\`ere un  rep\`ere $\{ u_1(x_1), \ldots, u_n(x_1) \}$ en $x_1$, de sorte que $u_1(x_1)$ appartienne \`a ${\mathcal F}_1(x_1)$, $\{ u_1(x_1), \ldots, u_{n-1}(x_1) \}$ soit une base de ${\mathcal F}_2(x_1)={\mathcal F}_1(x_1)^{\perp}$, et $u_n(x_1)$ soit isotrope avec de plus $g_{x_1}( u_{1}(x_1), u_{n}(x_1) )=1$.
 Les points $1(a)$ et $1(b)$ du Th\'eor\`eme \ref{thm.integrabilite-stable} fournissent alors une suite de   rep\`eres $((u_{1,k}(x_1), \ldots, u_{n,k}(x_1) ))$ en $x_1$, convergeant vers  $(u_1(x_1), \ldots, u_n(x_1) )$, telle que
\begin{eqnarray*}
|| D_{x_1} \phi_X^{s_k} (u_{1,k}(x_1)) || & = & \Theta(\mu_2(k)^2)  \\
|| D_{x_1} \phi^{s_k} (u_{i,k}(x_1)) || & = &\Theta( \mu_2(k)), \ i = 2, \ldots, n-1 \\
|| D_{x_1} \phi^{s_k} (u_{n,k}(x_1)) || & = & \Theta(1)
\end{eqnarray*}

\begin{lemme}
\label{lem.restriction.weyl}
On a les restrictions suivantes sur le tenseur de Weyl en $y$~:
\begin{enumerate}
\item{$\left. W \right|_{{\mathcal F}_2(x_1)}=0$.}
\item{Pour tout $u \in T_{x_1}M$, on a $W_{x_1}(u_1(x_1),u,u_1(x_1))=0$.}

\item{Pour $i,j,l \in \{1, \ldots, n  \}$, on a $W_{x_1}(u_i(x_1),u_j(x_1),u_l(x_1)) \in {\mathcal F}_1(x_1)$, sauf \'eventuellement si $i=l=n$ ou $j=l=n$.}
\end{enumerate}
\end{lemme}

\begin{proof}
 Comme,   $\phi_X^t.x_1$ reste dans un compact de $M$ pour $t \geq 0$, il existe $C > 0$ tel que
$$  || W_{\phi^t_X x_1} (a,b,c) || \leq C || a|| \cdot ||b|| \cdot || c || \ \forall  \ t \geq 0.$$

Si deux des indices $i,j,l$ sont dans $\{1, \ldots, n-1 \}$, alors~:
$$|| W_{\phi_X^{s_k}.x_1}(D_{x_1} \phi^{s_k}_X (u_{i,k}(x_1)), D_{x_1} \phi^{s_k}_X(u_{j,k}(x_1)), D_{x_1} \phi^{s_k}_X(u_{l,k}(x_1)) )|| = o\left( \mu_2(k)\right)$$
Alors par \'equivariance de $W$,
$$||D_{x_1} \phi^{s_k}_X( W_{x_1}(u_{i,k}(x_1), u_{j,k}(x_1), u_{l,k}(x_1))) ||= o\left(\mu_2(k) \right)$$
Par le Th\'eor\`eme \ref{thm.integrabilite-stable}, cela prouve que 

\begin{eqnarray}
\label{eqn.weyl.in.u1}
 W_{x_1}(u_i(x_1), u_j(x_1), u_l(x_1)) \in {\mathcal F}_1(x_1) 
\end{eqnarray}

De m\^eme, pour tout $u \in T_{x_1}M$~:
$$|| W_{\phi_X^{s_k}.x_1}(D_{x_1} \phi^{s_k}_X (u_{1,k}(x_1)), D_{x_1} \phi^{s_k}_X(u), D_{x_1} \phi^{s_k}_X(u_{1,k}(x_1)) )|| = o\left( \mu_2(k)^2 \right)$$ 
Le point $1(a)$ du Th\'eor\`eme \ref{thm.integrabilite-stable} permet alors d'affirmer que 
$$W_{x_1}(u_1(x_1),u,u_1(x_1))=0.$$

Enfin, le m\^eme type  d'argument montre que
\begin{eqnarray}
\label{eqn.weyl0.on.hyperplane}
W_{x_1}(u_i(x_1), u_j(x_1), u_l(x_1)) = 0  \qquad i,j,l < n.
\end{eqnarray} \end{proof}

Nous allons \`a pr\'esent montrer que le tenseur de Weyl s'annule en $x_1$.  Pour cela, on constate que d'apr\`es le choix des suites $(u_{i,k}(x_1))$, les vecteurs $\frac{1}{\mu_2(k)}D_{x_1} \phi_X^{s_k} (u_{i,k}(x_1))$ sont de normes born\'ees si $i<n$.  Quitte \`a consid\'erer une sous-suite de $(s_k)$, on aura que pour $i<n$, $\frac{1}{\mu_2(k)}D_{x_1} \phi_X^{s_k} (u_{i,k}(x_1))$ converge vers un vecteur $v_{i,\infty} \in T_{x_0}M$.  Par ailleurs, d'apr\`es la Proposition \ref{prop.limite.stable}, il existe une suite $(\nu_k)$ telle que $D_{x_1} \phi_X^{s_k} (\nu_k u_{1}(x_1))$ tende vers $u_{\infty} \in \BR^*.\Delta^{\prime}(0)$.

Par invariance du tenseur de Weyl, et par le point $(2)$ du Lemme \ref{lem.restriction.weyl}, on a pour tout $k \in \BN$~:
$$ W_{\phi^{s_k}_X.x_1}( D_{x_1} \phi^{s_k}_X(\nu_k u_1(x_1)), \frac{1}{\mu_2(k)}  D_{x_1} \phi^{s_k}_X(u_{i,k}(x_1)),  D_{x_1} \phi^{s_k}_X(\nu_k u_1(x_1)))=0.$$
En passant \`a la limite, on obtient que~:
\begin{eqnarray}
\label{eqn.weyl0}
W_{x_0}(u_{\infty}, v_{i,\infty}, u_{\infty})=0
\end{eqnarray}

Maintenant, par le Th\'eor\`eme \ref{thm.integrabilite-stable}, et par la Proposition \ref{prop.limite.stable}, il existe $\delta>0$ tel que $D_{x_1}\phi_X^{s_k}(u_n(x_1)) \to \delta u_{\infty}$.  On obtient alors~:

\begin{eqnarray*}
& & \frac{1}{\mu_2(k)} D_{x_1} \phi^{s_k}_X (W_{x_1} (u_n(x_1), u_{i,k}(x_1),u_n(x_1))) \\
& = & W_{\phi^{s_k}_X.x_1}( D_{x_1} \phi^{s_k}_X(u_n(x_1)), \frac{1}{\mu_2(k)}  D_{x_1} \phi^{s_k}_X(u_{i,k}(x_1)),  D_{x_1} \phi^{s_k}_X(u_n(x_1))) \\
& \rightarrow & \delta^2 W_{x_0}(u_\infty, v_{i,\infty}, u_\infty) \\
& = & 0
\end{eqnarray*}
par l'\'equation (\ref{eqn.weyl0}).

La suite des vecteurs  $(W_{x_1}(u_n(x_1), u_{i,k}(x_1), u_n(x_1)))$ est contract\'ee  par $(D_{x_1} \phi^{s_k}_X)$ plus vite que $(\mu_2(k))$.  Donc, par le Th\'eor\`eme \ref{thm.integrabilite-stable}, on a~:
$$ W_{x_1}(u_n(x_1), u_i(x_1), u_n(x_1)) \in {\mathcal F}_1(x_1).$$

Cette relation, jointe au Lemme \ref{lem.restriction.weyl}, permet de  conclure que l'image de $W_{x_1}$ est contenue dans $\BR. u_1(x_1)={\mathcal F}_1(x_1)$ et que le tenseur de Weyl est nul sur $u_1(x_1)^\perp={\mathcal F}_2(x_1)$.  Comme on a par ailleurs~:
\begin{eqnarray*}
 g_{x_1}(W_{x_1}(u_i(x_1), u_j(x_1),u_n(x_1)), u_n(x_1) ) & =& \\
  - g_{x_1}( W_{x_1}(u_i(x_1), u_j(x_1), u_n(x_1)), u_n(x_1) ) &=& 0,
  \end{eqnarray*} 

on obtient que  la composante de Weyl sur $u_1(x_1)$ est aussi nulle. La conclusion est que le tenseur de  Weyl s'annule en $x_1$, et aussi en $x_0$ puisque $x_0=\lim_{k \to \infty}\phi_X^{s_k}.x_1$.

\subsubsection{Quatri\`eme \'etape : annulation du tenseur de Weyl sur $\Delta$ et  fin de la preuve}

Consid\'erons \`a nouveau un param\'etrage de $\Lambda$ au voisinage de $o$ de la forme  $\Lambda(s):= \pi_G(e^{s \xi_1})$, o\`u $\xi_1 \in \lien^{-}$.  Comme $\{h^t \}$ fixe tous les points de $\Lambda$, il en va de m\^eme de $\{ h_v^t \}:=e^{-v \xi_1}h^te^{v \xi_1}$,  pour tout $v \in \BR$.  En particulier $\{ h_v^t \} <P$, et la relation $h^te^{v\xi_1}h_v^{-t}=e^{v \xi_1}$ montre que $\{ h_v^t \}$ est un flot d'holonomie de $\{ h^t \}$ en $\Lambda(v)$.  Pour  $v \in (-\epsilon, \epsilon)$ avec $\epsilon>0$ assez petit, $v \mapsto \pi(\exp(b_0,v \xi_1))$ est un param\'etrage de $\Delta$ au voisinage de $x_0$.  Il n'est pas tr\`es difficile de v\'erifier que $\{ h_v^t \}$ est un flot d'holonomie de $\{ \phi_X^t \} $ en $\Delta(v)$ (voir par exemple la Proposition 4.3 de \cite{nilpotent}).    Les lemmes et propositions montr\'es lors des trois premi\`eres \'etapes restent alors valables si l'on remplace $o$ par $\Lambda(v)$, $x_0$ par $\Delta(v)$,  et la g\'eod\'esique $\beta$ par $e^{-v\xi_1}.\beta$.  
En particulier, on obtient de la m\^eme mani\`ere que pr\'ec\'edemment que le tenseur de Weyl est nul en $\Delta(v)$ pour tout $v \in (-\epsilon, \epsilon)$, en ayant choisi un $\epsilon$ assez petit.

Soit $V$ le voisinage original de $x_1 = \alpha(1)$.  La suite $( \phi^{s_k}_X )$ est stable sur $V$, et les images ($\phi^{s_k}_X. V)$ tendent vers un segment de $\Delta$, o\`u le tenseur de Weyl s'annule.  On conclut par \cite[Prop 4(i)]{frances-causal} que $W$ s'annule sur $V$.  Comme $M$ est suppos\'ee analytique, alors elle est conform\'ement plate partout.

\section{Un flot contractant sur une vari\'et\'e non conform\'ement plate}
\label{sec.exemple-contraction}

Le fait qu'un flot conforme contracte un ouvert non vide assure la platitude conforme de l'ouvert lorsque la signature est riemannienne ou lorentzienne, mais ne suffit pas forc\'ement  pour prouver l'annulation du tenseur de Weyl sur cet ouvert en signature arbitraire.  Comme on l'a vu plus haut, les   vitesses relatives de contractions  dans les diff\'erentes directions jouent un r\^ole important.  On peut n\'eanmoins \'enoncer le r\'esultat g\'en\'eral suivant~:

\begin{theor}\cite[Th. 1.2]{frances-degenere}
Soit $(M,g)$ une vari\'et\'e pseudo-riemannienne de dimension $n \geq 3$.  On suppose qu'il existe un ouvert non vide $U \subset M$ et une suite de plongements conformes $f_k:U \to M$ telle que la suite d'ouverts $f_k(U)$ converge vers un point de $M$ pour la topologie de Hausdorff. Alors $(U,g)$ est localement conform\'ement Ricci-plat. 
\end{theor}
Par localement conform\'ement Ricci-plat, on entend qu'il existe au voisinage de chaque point de $U$ une m\'etrique Ricci-plate dans la classe conforme de $g$.

Dans l'exemple qui suit, nous construisons un flot  lin\'eaire  qui contracte un voisinage de l'origine, et agit conform\'ement pour une  m\'etrique {\it qui n'est pas conform\'ement plate}.  

Consid\'erons sur $\BR^6$ la m\'etrique
$$ g= \dd x_1 \dd x_6 + \dd x_2 \dd x_5 + \dd x_3 \dd x_4 + x_1 x_2 x_3 \dd x_1^2$$

Le flot lin\'eaire
$$
\left( \begin{array}{cccccc}
e ^{-t} &  &  &  &  &   \\
    & e^{-2t} &  &  &  &  \\
  &   &  e^{-3t} &  &  &  \\
  &   &    &  e^{-5t}  &  &  \\
  &   &  &   & e^{-6t}  &    \\
  &   &  &   &   &   e^{-7t}
\end{array} \right)
$$

est conforme pour cette m\'etrique.

On va montrer que la courbure de Weyl de $g$ ne s'annule pas sur le triplet $(e_1, e_2, e_1)$, o\`u les $e_i$ forment la base standard de $\BR^6$.  Par abus de langage, on notera  encore $e_i$ le champ de vecteurs invariant par translations, valant $e_i$ en $0$.  On a  $[e_i,e_j] = 0$ pour tous $i,j$.  On calcule alors
\begin{eqnarray*}
\langle \nabla_{e_1} e_1, e_1 \rangle & = & \frac{1}{2} \frac{\partial}{\partial x_1} \langle e_1, e_1 \rangle = \frac{1}{2} x_2 x_3 \\
\langle \nabla_{e_1} e_1, e_2 \rangle & = & - \frac{1}{2} \frac{\partial}{\partial x_2} \langle e_1, e_1 \rangle = - \frac{1}{2} x_1 x_3 \\
\langle \nabla_{e_1} e_1 , e_3 \rangle & = & - \frac{1}{2} \frac{\partial}{\partial x_3} \langle e_1, e_1 \rangle = - \frac{1}{2} x_1 x_2 \\
\langle \nabla_{e_1} e_1, e_i \rangle & = & 0 \qquad i = 4, 5, 6 
\end{eqnarray*}

Donc 
$$\nabla_{e_1} e_1 = - x_1 x_2 e_4 - x_1 x_3 e_5 + x_2 x_3 e_6$$

Des calculs similaires donnent
\begin{eqnarray*}
\nabla_{e_1} e_2 & = & - x_1 x_3 e_6 \\
\nabla_{e_1} e_3 & = & - x_1 x_2 e_6 \\
\nabla_{e_j} e_i & = & 0 \qquad j = 1, \ldots, 6, \ i = 4, 5, 6
\end{eqnarray*}

Puis
\begin{eqnarray*}
\nabla_{e_2} \nabla_{e_1} e_1 & = & - x_1 e_4 + x_3 e_6 \\
\nabla_{e_1} \nabla_{e_1} e_2 & = & - x_3 e_6
\end{eqnarray*} 

Le tenseur de courbure $R$ de $g$ prend alors la valeur
\begin{eqnarray*}
R(e_1, e_2)e_1 & = & (\nabla_{e_2} \nabla_{e_1} - \nabla_{e_1} \nabla_{e_2}) e_1 \\
& = & (\nabla_{e_2} \nabla_{e_1} - \nabla_{e_1} \nabla_{e_1}) e_2 \\
& = & - x_1 e_4 + 2 x_3 e_6  
\end{eqnarray*}

Supposons que $x_2 x_3 = 0$ mais $x_1 \neq 0$.  Alors
$$\langle W(e_1,e_2)e_1, e_3 \rangle = \langle R(e_1, e_2)e_1, e_3 \rangle - (g \circ P)(e_1,e_2, e_1, e_3)$$

o\`u $\circ$ est le produit de Kulkarni-Nomizu, et $P$ est le tenseur de Schouten (voir \cite[p.48 ]{besse}).  Ce produit est
\begin{eqnarray*}
 g \circ P (e_1, e_2, e_1, e_3) & = & g(e_1,e_1)P(e_2, e_3)+ g(e_2,e_3)P(e_1,e_1)  \\
&  &    - g(e_1,e_3)P(e_2,e_1) - g(e_2,e_1)P(e_1,e_3) \\
& = & x_1x_2x_3 P(e_2,e_3) \\
& = & 0
\end{eqnarray*}

Donc
$$ \langle W(e_1,e_2)e_1,e_3 \rangle = \langle R(e_1,e_2)e_1,e_3 \rangle = - x_1 \neq 0 $$
et on conclut que $W(e_1, e_2, e_1)$ ne s'annule sur aucun voisinage de l'origine.



\begin{thebibliography}{00}

\bibitem[A]{aleks}  D.V Alekseevskii~: Groups of conformal transformations of Riemannian spaces. Math. USSR, Sb. 18 (1972), 285-301. 

\bibitem[BCDGM]{primer} T. Barbot, V. Charette, T. Drumm, W. M. Goldman, K. Melnick :
  A primer on the $(2+1)$ Einstein universe, in {Recent developments in pseudo-Riemannian geometry} (D. Alexeevskii and H. Baum, eds). Z\"urich: European Mathematical Society. ESI Lectures in Mathematics and Physics, 179-229 (2008).

\bibitem[BCH]{beig} M.S Capocci, R. Beig, G. S. Hall :  Zeros of conformal vector fields, {  Classical Quantum Gravity}  { 14}  (1997),  no. 3, 49--52.

\bibitem[B]{besse}{ A.Besse} ~:  {Einstein manifolds}. Ergebnisse der Mathematik und ihrer Grenzgebiete (3), 10. Springer-Verlag, Berlin, 1987.

\bibitem[\v{C}M]{cap.me.parabolictrans} A. {\v C}ap and K. Melnick: Essential Killing fields of parabolic geometries, arXiv:1208.5510 (2012).

\bibitem[Ca1]{capocci1} M. S. Capocci :  Essential conformal vector fields, {  Classical Quantum Gravity}  { 16}  (1999),  no. 3, 927--935.

\bibitem[Ca2]{capocci2} M. S. Capocci :  Conformal vector fields and non-degenerate distributions, {  Classical Quantum Gravity}  { 13}  (1996), 1717--1726.


\bibitem[F1]{ferrand} J. Ferrand : The action of conformal transformations on a Riemannian manifold, { Math. Ann.}  { 304}  (1996),  no. 2, 277--291. 

\bibitem[F2]{ferrand2} J. Ferrand~: Sur un lemme d'Alekseevskii relatif aux transformations conformes.  C.  R. Acad. Sci., Paris, S\'er. A 284, 121-123 (1977).


\bibitem[Fr1]{frances-causal} C. Frances :  Causal conformal vector fields, and singularities of twistor spinors, {  Ann. Global Anal. Geom.}  { 32}  (2007),  no. 3, 277--295. 

\bibitem[Fr2]{charlesthese} C. Frances : thesis. http://www.math.u-psud.fr/~frances/


\bibitem[Fr3]{frances-riemannien} C. Frances :  Local dynamics of conformal vector fields, {  Geometriae Dedicata}. { 158} (2012), no. 1, 35--59.

\bibitem[Fr4]{frances-degenere} C. Frances~: D\'eg\'enerescence locale des transformations pseudo-riemanniennes conformes. arXiv:1008.2436v1.  \`A para\^{\i}tre dans  Annales de l'Institut Fourier. 



\bibitem[FM]{nilpotent} C. Frances, K. Melnick :  Conformal actions of nilpotent groups on pseudo-Riemannian manifolds, {Duke Math. J.}  {153} (2010), no. 3, 511--550.

\bibitem[G]{gromov.rgs} M. Gromov~: Rigid transformations groups, in {Geometrie Differentielle} (D. Bernard and Y. Choquet-Bruhat, eds.) Paris: Hermann, 1988.
   
\bibitem[Ko]{kobayashi} S. Kobayashi~: Transformation groups in differential geometry.  Berlin: Springer-Verlag (1995).

   
\bibitem[KR1]{kuehnel} W. K\"uhnel, H. B. Rademacher :  Essential conformal fields in pseudo-Riemannian geometry. II, {  J. Math. Sci. Univ. Tokyo} {  4}  (1997),  no. 3, 649--662.


\bibitem[KR2]{kuehnel2} W. K\"uhnel, H. B. Rademacher : Essential conformal fields in pseudo-Riemannian geometry, {J. Math. Pures Appl.} (9)  {74}  (1995),  no. 5, 453--481. 

\bibitem[KR3]{kuehnel4} W. K\"uhnel, H.B. Rademacher : Conformal vector fields on pseudo-Riemannian spaces. {Differential Geom. Appl.}  {7} (1997), no. 3, 237--250.


\bibitem[Me]{me.frobenius} K. Melnick~: A Frobenius theorem for Cartan geometries, with applications. 	\emph{eprint arXiv:0812.0624v2}. L'Enseignement Math\'ematique (S\'er. II) 57 no. 1-2 (2011) 57-89.

\bibitem[Mo]{morris} D. W. Morris~: Introduction to arithmetic groups, \emph{eprint arXiv:math/0106063v3}.

\bibitem[NO]{nagano} T, Nagano, T. Ochiai~:
On compact Riemannian manifolds admitting essential projective transformations.  J. Fac. Sci., Univ. Tokyo, Sect. I A 33, 233-246 (1986).


\bibitem[Sh]{sharpe}  R. W. Sharpe : {Differential Geometry : Cartan's
    generalization of Klein's Erlangen Program}. Berlin: Springer-Verlag (1997).

\bibitem[St]{steller} M. Steller :  Conformal vector fields on spacetimes, {  Ann. Global Anal. Geom.}  { 29}  (2006),  no. 4, 293--317.

\bibitem[Yo76]{yoshi76} Y. Yoshimatsu~: On a theorem of Alekseevskii concerning conformal transformations. J. Math. Soc. Japan {28} (1976), no. 2, 278--289.
\end{thebibliography}
\end{document}